\documentclass{amsart}

\usepackage[latin1]{inputenc}
\usepackage{eepic,epic, amsmath, mathtools}
\usepackage{graphicx}
\usepackage{xypic}

\newcommand{\n}{*=0{{\scriptstyle\bullet}}}
\newcommand{\eru}{\ar@{-}[ru]} 
\newcommand{\erd}{\ar@{-}[rd]}
\newcommand{\er}{\ar@{-}[r]} 
\newcommand{\eu}{\ar@{-}[u]} 

\newcommand{\sym}{\mathcal{S}} 

\newcommand{\NN}{{\mathbb N}}  
\newcommand{\dyck}{{\mathcal D}} 

\DeclareMathOperator{\Knuth}{\mathrm{Knuth}}
\DeclareMathOperator{\Ri}{\mathrm{Richards}}
\DeclareMathOperator{\KRi}{\mathrm{Knuth-Richards}}
\DeclareMathOperator{\KRo}{\mathrm{Knuth-Rotem}}
\DeclareMathOperator{\Astrid}{\mathrm{Reifegerste}}
\DeclareMathOperator{\MDD}{\mathrm{Mansour-Deng-Du}}
\DeclareMathOperator{\SiS}{\mathrm{Simion-Schmidt}}
\DeclareMathOperator{\Krattan}{\mathrm{Krattenthaler}}
\DeclareMathOperator{\ED}{\mathrm{Elizalde-Deutsch}}
\DeclareMathOperator{\West}{\mathrm{West}}

\DeclareMathOperator{\reverse}{\mathrm{reverse}} 
 
\DeclareMathOperator{\inverse}{\mathrm{inverse}} 
\DeclareMathOperator{\identity}{\mathrm{identity}} 

\DeclareMathOperator{\comp}{\mathrm{comp}}
\DeclareMathOperator{\compR}{\mathrm{comp.\!r}}

\DeclareMathOperator{\rank}{\mathrm{rank}}
\DeclareMathOperator{\rankR}{\mathrm{rank.\!r}}
\DeclareMathOperator{\rankRC}{\mathrm{rank.\!r.\!c}}
\DeclareMathOperator{\Nrank}{\mathrm{n-rank}}
\DeclareMathOperator{\lmax}{\mathrm{lmax}}
\DeclareMathOperator{\rmax}{\mathrm{rmax}}
\DeclareMathOperator{\lmin}{\mathrm{lmin}}
\DeclareMathOperator{\rmin}{\mathrm{rmin}}
\DeclareMathOperator{\Nrmin}{\mathrm{n-rmin}}
\DeclareMathOperator{\ldr}{\mathrm{ldr}}
\DeclareMathOperator{\ldrR}{\mathrm{ldr.\!r}}
\DeclareMathOperator{\ldrI}{\mathrm{ldr.\!i}}

\DeclareMathOperator{\Mldr}{\mathrm{m-ldr}}
\DeclareMathOperator{\MldrI}{\mathrm{m-ldr.\!i}}
\DeclareMathOperator{\lir}{\mathrm{lir}}
\DeclareMathOperator{\lirR}{\mathrm{lir.\!r}}
\DeclareMathOperator{\lirI}{\mathrm{lir.\!i}}

\DeclareMathOperator{\rdr}{\mathrm{rdr}}
\DeclareMathOperator{\rdrI}{\mathrm{rdr}.\!i}
\DeclareMathOperator{\rir}{\mathrm{rir}}
\DeclareMathOperator{\rirI}{\mathrm{rir.\!i}}
\DeclareMathOperator{\fp}{\mathrm{fix}}
\DeclareMathOperator{\cyc}{\mathrm{cyc}}
\DeclareMathOperator{\zeil}{\mathrm{zeil}}
\DeclareMathOperator{\zeilC}{\mathrm{zeil.\!c}}
\DeclareMathOperator{\head}{\mathrm{head}}
\DeclareMathOperator{\headI}{\mathrm{head.\!i}}
\DeclareMathOperator{\headIR}{\mathrm{head.\!i.\!r}}
\DeclareMathOperator{\last}{\mathrm{last}}
\DeclareMathOperator{\lastI}{\mathrm{last.\!i}}
\DeclareMathOperator{\lastIR}{\mathrm{last.\!i.\!r}}
\DeclareMathOperator{\lis}{\mathrm{lis}}
\DeclareMathOperator{\lds}{\mathrm{lds}}

\DeclareMathOperator{\exc}{\mathrm{exc}}
\DeclareMathOperator{\excR}{\mathrm{exc.\!r}}
\DeclareMathOperator{\peak}{\mathrm{peak}}
\DeclareMathOperator{\peakI}{\mathrm{peak.\!i}}

\DeclareMathOperator{\peakR}{\mathrm{peak.\!r}}
\DeclareMathOperator{\valley}{\mathrm{valley}}
\DeclareMathOperator{\valleyC}{\mathrm{valley.\!c}}
\DeclareMathOperator{\valleyI}{\mathrm{valley.\!i}}

\DeclareMathOperator{\des}{\mathrm{des}} 
\DeclareMathOperator{\Ndes}{\mathrm{n-des}} 
\DeclareMathOperator{\asc}{\mathrm{asc}}
\DeclareMathOperator{\ascR}{\mathrm{asc.\!r}}
\DeclareMathOperator{\slmax}{\mathrm{slmax}}
\DeclareMathOperator{\slmaxI}{\mathrm{slmax.\!i}}
\DeclareMathOperator{\slmaxIR}{\mathrm{slmax.\!i.\!r}}
\DeclareMathOperator{\slmaxRI}{\mathrm{slmax.\!r.\!i}}
\DeclareMathOperator{\slmaxR}{\mathrm{slmax.\!r}}
\DeclareMathOperator{\slmaxC}{\mathrm{slmax.\!c}}
\DeclareMathOperator{\slmaxRC}{\mathrm{slmax.\!r.\!c}}
\DeclareMathOperator{\slmaxRCI}{\mathrm{slmax.\!r.\!c.\!i}}
\DeclareMathOperator{\slmaxCR}{\mathrm{slmax.\!c.\!r}}

\DeclareMathOperator{\LMIN}{\mathfrak{lmin}}
\DeclareMathOperator{\RMIN}{\mathfrak{rmin}}
\DeclareMathOperator{\LMAX}{\mathfrak{lmax}}
\DeclareMathOperator{\RMAX}{\mathfrak{rmax}}
\DeclareMathOperator{\PEAK}{\mathfrak{peak}}

\newcommand{\CK}{\ensuremath{\Phi}}
\newcommand{\STAT}{\ensuremath{\mathrm{STAT}}}
\newcommand{\emptyword}{\epsilon}

\newtheorem{theorem}{Theorem}

\newtheorem{lemma}[theorem]{Lemma}

\newtheorem{definition}[theorem]{Definition}

\newcommand\p{\circle*{0.25}}

\author{Anders Claesson\and Sergey Kitaev} 
\title[Classification of bijections] 
      {Classification of bijections between \\ 321- and 132-avoiding permutations}
\thanks{The research presented here was supported by grant no. 060005012 from
the Icelandic Research Fund}
\address{The Mathematics Institute, Reykjavik University, Kringlan 1,
103 Reykjavik, Iceland} 

\begin{document}
\maketitle

\begin{abstract}
  It is well-known, and was first established by Knuth in 1969, that
  the number of $321$-avoiding permutations is equal to that of
  $132$-avoiding permutations. In the literature one can find many
  subsequent bijective proofs of this fact. It turns out that some of
  the published bijections can easily be obtained from others. In this
  paper we describe all bijections we were able to find in the
  literature and show how they are related to each other via
  ``trivial'' bijections.  We classify the bijections according to
  statistics preserved (from a fixed, but large, set of statistics),
  obtaining substantial extensions of known results. Thus, we give a
  comprehensive survey and a systematic analysis of these bijections.
  
  We also give a recursive description of the algorithmic bijection
  given by Richards in 1988 (combined with a bijection by Knuth from
  1969). This bijection is equivalent to the celebrated bijection of
  Simion and Schmidt (1985), as well as to the bijection given by
  Krattenthaler in 2001, and it respects 11 statistics---the largest
  number of statistics any of the bijections respects.
\end{abstract}

\setcounter{tocdepth}{1}
\tableofcontents

\section{Introduction and main results}
\label{sec:in}

Given two different bijections between two sets of combinatorial
objects, what does it mean to say that one bijection is better than
the other? Perhaps, a reasonable answer would be ``The one that is
easier to describe.'' While the ease of description and how easy it is
to prove properties of the bijection using the description is one
aspect to consider, an even more important aspect, in our opinion, is
how well the bijection reflects and translates properties of elements
of the respective sets.

A natural measure for a bijection between two sets of permutations,
then, is how many statistics the bijection preserves. Obviously, we
don't have an exhaustive list of permutation statistics, but we have
used the following list as our ``base'' set:
\begin{gather*}
\asc, \des, \exc, \ldr, \rdr, \lir, \rir, \zeil, \comp, \lmax, \lmin,
\rmax, \rmin,\\
\head, \last, \peak, \valley, \lds, \lis, \rank, \cyc, \fp, \slmax.
\end{gather*}
These are defined in Section~\ref{sec:stat}. To make sure we find all
statistics that a given bijection ``essentially'' preserves, we
include in our list of statistics those that are obtained from our
``basic'' statistics by applying to them the {\em trivial bijections}
on permutations ({\em reverse=r}, {\em complement=c}, {\em
  inverse=i\/}) and their compositions. Moreover, for each statistic
stat, in this extended list we consider two other statistics:
$\textrm{n-stat}(\pi) = n-\textrm{stat}(\pi)$ and
$\textrm{m-stat}(\pi) = n+1-\textrm{stat}(\pi)$, where $n$ is the
length of the permutation. The meaning of n-stat or m-stat is
often ``non-stat''; for example, n-fix counts non-fixed-points.

This way each basic statistic gives rise to 24 statistics. The base
set contains 23 statistics, giving a total of 552 statistics. There
are, however, many statistics in that set that are equal as functions;
for instance, $\des = \ascR$, and $\peak = \peakR = \valleyC$, where
we use a dot to denote composition of functions. Choosing one
representative from each of the classes of equal statistics results in
a final set of $190$ statistics; we call this set \STAT.  In practice
we settled for ``empirical equality'' when putting together \STAT: we
considered two statistics equal if they gave the same value on all
5914 permutations of length at most 7.


In the theorems below, the statistics presented are linearly
independent. An example of linear dependence among the statistics
over permutations avoiding $132$ is $\lmin -\lmax + \Ndes -\head = 0$.
The results below are also maximal in that they cannot be
non-trivially extended using statistics from $\STAT$. That is, adding
one more pair of equidistributed statistics from \STAT\ to any of the
results would create a linear dependency among the statistics.

A permutation $\pi=a_1a_2\dots a_n$ avoids the \emph{pattern} 321 if
there are no indices $i<j<k$ such that $a_k<a_j<a_i$. It avoids 132 if
there are no indices $i<j<k$ such that $a_i<a_k<a_j$. Avoidance of
other patterns is defined similarly.

Knuth~\cite{K69,K73} showed that the number of permutations avoiding a
pattern of length 3 is independent of the pattern.  This number is the
{\em $n$-th Catalan number}, $C_n=\frac{1}{n+1}{2n \choose n}$.  To
prove this it suffices, due to the symmetry afforded by the trivial
bijections on permutations, to consider one representative from
$\{123,321\}$ and one from $\{132,231,213,312\}$. That symmetry also
means that to prove this bijectively, it suffices to find a bijection
from the set of permutations avoiding a pattern in one of the classes
to permutations avoiding a pattern in the other.  This turns out to be
a rather complicated problem. Several authors have, however, succeeded
in doing so \cite{ED03,EP04,K01,R02,R88,R75,SS85,W95}; we call those
bijections
\begin{gather*}
  \Knuth,\, \KRo,\, \SiS,\, \KRi,\, \West,\, \\ \Krattan,\,
  \Astrid,\, \ED,\, \text{ and }\MDD.
\end{gather*}
They are described in Section~\ref{sec:bij}. In Section~\ref{sec:CK}
we define, using recursion, a ``new'' bijection called $\CK$. It turns
out that $\CK$ is related, via trivial bijections, to the bijection by
Knuth and Richards.

The main results of this paper are contained in the following three
theorems. The first theorem substantially extends what was previously
known about statistics preserved by the bijections. In bold we mark
the results known before this paper---there is a total of 7 pairs of
those; the remaining 68 pairs are new.

\begin{theorem}\label{extensions}
  The following results are maximal in the sense that adding one more
  pair of equidistributed statistics from \STAT\ to any of the results
  would create a linear dependency among the statistics. In bold we
  mark the results that were known; also, we indicate the sets between
  which a bijection acts.\medskip

  {\smaller
  \begin{list}{}{
      \setlength{\leftmargin}{1ex}
      \setlength{\itemsep}{2.5ex}
    }

  \item[$(11)$] {\bf Knuth-Richards}, 132-avoiding permutations $\rightarrow$ 123-avoiding permutations \medskip\\
    $\begin{array}{lllllllllll}
\valleyI & \valley  & \lmin & \ldrI & \headI & \compR & \rank & \ldr  & \lirI   & \lir     & \rmax   \\
\valley  & \valleyI & \lmin & \ldr  & \head  & \compR & \rank & \ldrI & \slmaxC & \slmaxIR & \headIR  
    \end{array}$
  \item[$(11)$] {\bf Simion-Schmidt}, 123-avoiding permutations
    $\rightarrow$ 132-avoiding permutations\medskip\\
    $\begin{array}{lllllllllll}
\valley & \valleyI & {\bf lmin} & \ldr & \head & \compR & \rank & \ldrI & \slmaxC & \slmaxIR & \headIR \\ 
\valley & \valleyI & {\bf lmin} & \ldr & \head & \compR & \rank & \ldrI & \lir    & \lirI     & \rmin 
  \end{array}$
  \item[$(11)$] {\bf Krattenthaler}, 123-avoiding permutations
    $\rightarrow$ 132-avoiding permutations\medskip\\
    $\begin{array}{lllllllllll}
 \peakI & \peak     & \rmax & \zeil & \lastIR & \compR & \rankRC & \rdr  & \slmaxRI & \slmaxR & \last \\
 \valley & \valleyI & \lmin & \ldr  & \head   & \compR & \rank   & \ldrI & \lir     & \lirI   & \rmin 
    \end{array}$
 \item[$(11)$] {\bf Mansour-Deng-Du}, 321-avoiding permutations
    $\rightarrow$ 231-avoiding permutations\medskip\\
    $\begin{array}{lllllllllll}
   \valley & \peakI & \rmin & \rir & \last & \comp & \rankR & \lirI & \slmaxCR & \slmaxI & \headI \\
   \valley & \peakI & \rmin & \rir & \last & \comp & \rankR & \lirI & \rdr     & \ldrI   & \lmin
  \end{array}$
\item[$(9)$] {\bf Knuth-Rotem}, 321-avoiding permutations $\rightarrow$
  132-avoiding permutations\medskip\\
  $\begin{array}{lllllllll}
    \valleyI & \peak   & {\bf exc} & \slmax & \head & \slmaxRCI & \rirI & \lir  & \lastI \\ 
    \valleyI & \valley & {\bf des} & \rdr   & \ldrI & \zeil     & \lmax & \rmin & \Mldr
  \end{array}$
\item[$(9)$] {\bf Reifegerste}, 321-avoiding permutations $\rightarrow$
  132-avoiding permutations\medskip\\
  $\begin{array}{lllllllll}
    \valley & \peakI & {\bf exc} & \slmaxI & \headI & \slmaxRC & \rir & \lirI & \last \\ 
    \valley & \valleyI & {\bf des} & \zeil &\ldr & \rdr & \rmin & \lmax & \MldrI
  \end{array}$
\item[$(7)$] {\bf West}, 123-avoiding permutations $\rightarrow$
  132-avoiding permutations\medskip\\
  $\begin{array}{lllllll}
    \valleyI & \excR & \slmaxIR & \slmaxC & \ldr  & \ldrI & \head \\
    \valleyI & \asc  & \lirI    & \comp   & \rmax & \ldrI & \head
  \end{array}$
\item[$(5)$] {\bf Knuth}, 321-avoiding permutations $\rightarrow$
  132-avoiding permutations\medskip\\
  $\begin{array}{lllll}
    {\bf exc} & {\bf fix} & \lirI & \lir & {\bf lis}\\ 
    {\bf exc} & {\bf fix} & \rmin & \lmax & \text{\bf n-rank}
  \end{array}$
\item[$(1)$] {\bf Elizalde-Deutsch}, 321-avoiding permutations
  $\rightarrow$ 132-avoiding permutations\medskip\\
  $\begin{array}{l}
    {\bf fix} \\
    {\bf fix} 
  \end{array}$
  \end{list}
  }
\end{theorem}

The numbers in parenthesis in Theorem~\ref{extensions} indicate the
number of statistics respected. It turns out that bijections with the
same number are ``trivially'' related. The next theorem makes this
precise.

\begin{theorem}\label{relations}
  The following relations among bijections between $321$- and
  $132$-avoiding permutations hold:
  \begin{align*}
    \reverse\,\circ\,\CK^{-1}\ 
    &=\ \inverse\,\circ\,\SiS\,\circ\,\reverse \\
    &=\ \inverse\,\circ\,\Krattan\,\circ\,\reverse\,\circ\,\inverse
    \\
    &=\ \inverse\,\circ\,\reverse\,\circ\,\MDD
    \\
    &=\ \KRi^{-1} \circ\,\reverse
    \shortintertext{and}
    \Astrid\ &=\ \inverse \,\circ\,\KRo\,\circ\,\inverse
  \end{align*}
  Also, there are no other relations among the bijections and their
  inverses via the trivial bijections that does not follow from the
  ones above.
\end{theorem}

Thus if we regard all bijections as bijections from $321$- to
$132$-avoiding permutations---applying the transformations in
Theorem~\ref{relations}---then we get the following condensed version
of Theorem~\ref{extensions}.

\begin{theorem}\label{statistics}
For bijections from 321- to 132-avoiding permutations we have the
following equidistribution results. These results are maximal in the
sense that adding one more pair of equidistributed statistics from
\STAT\ to any of the results would create a linear dependency among
the statistics.

{\smaller
  \begin{list}{}{
      \setlength{\leftmargin}{1ex}
      \addtolength{\itemsep}{2ex}
    }
  \item[$(11)$]{\bf \CK, Knuth-Richards, Krattenthaler, Mansour-Deng-Du,
    Simion-Schmidt}\medskip\\
    $\begin{array}{lllllllllll} 
\valley & \peakI   & \rmin & \rir & \last & \comp  & \rankR & \lirI & \slmaxCR & \slmaxI & \headI \\
\valley & \valleyI & \lmin & \ldr & \head & \compR & \rank  & \ldrI & \lir     & \lirI   & \rmin
    \end{array}$
  \item[$(9)$] {\bf Knuth-Rotem, Reifegerste}\medskip\\
    $\begin{array}{lllllllll} 
      \valley & \peakI & \exc & \slmaxI & \headI & \slmaxRC & \rir & \lirI & \last \\
      \valley & \valleyI & \des & \zeil & \ldr & \rdr & \rmin & \lmax &
      \MldrI
    \end{array}$
  \item[$(7)$] {\bf West}\medskip\\
    $\begin{array}{lllllll} 
      \peakI & \exc & \slmaxI & \slmaxRC & \rir & \lirI & \last \\ 
      \valleyI & \asc & \lirI & \comp & \rmax & \ldrI\ & \head
    \end{array}$
    \item[$(5)$] {\bf Knuth}\medskip\\
      $\begin{array}{lllll} 
        \exc & \fp & \lirI & \lir & \lis \\ 
        \exc & \fp & \rmin & \lmax & \Nrank
      \end{array}$
    \item[$(1)$] {\bf Elizalde-Deutsch} \medskip\\
      $\begin{array}{l} 
        \fp \\
        \fp 
      \end{array}$
  \end{list}
}
\end{theorem}

In Section~\ref{sec:stat} we define the relevant statistics; in
Section~\ref{sec:bij} we describe the bijections; in
Section~\ref{sec:CK} we give a new recursive description of the
bijection by Knuth and Richards; in Section~\ref{sec:relations} we
prove Theorem~\ref{relations}; and in Section~\ref{sec:statistics} we
prove Theorem~\ref{extensions}.

\section{Permutation statistics}
\label{sec:stat}
The permutation $\pi$ on $\{1,2,\dots,n\}$ that sends $1$ to
$a_1$, $2$ to $a_2$, etc, we denote $\pi=a_1a_2\dots a_n$, and we call
$a_i$ the $i$-th letter of $\pi$. A permutation statistic is simply a
function from permutations to $\NN$. For example, the permutation
statistic $\asc$ is defined thus: An \emph{ascent} in $\pi$ is a
letter that is followed by a larger letter; in other words, an $a_i$
such that $a_i<a_{i+1}$. By $\asc(\pi)$ we denote the number of
ascents in $\pi$. Similarly, a \emph{descent} is a letter followed by a
smaller letter, and by $\des(\pi)$ we denote the number of descents in
$\pi$.

For words $\alpha$ and $\beta$ over the alphabet $\NN$ we define that
$\alpha\prec\beta$ if for all letters $a$ in $\alpha$ and all letters
$b$ in $\beta$ we have $a < b$. For instance, $42\prec 569$.
A \emph{component} of $\pi$ is a nonempty segment $\tau$ of $\pi$ such
that $\pi = \sigma\tau\rho$ with $\sigma\prec\tau\prec\rho$, and such
that if $\tau$ = $\alpha\beta$ and $\alpha\prec\beta$ then $\alpha$ or
$\beta$ is empty. By $\comp(\pi)$ we denote the number of components
of $\pi$. For instance, $\comp(213645)=3$, the components being $21$,
$3$, and $645$.

A \emph{left-to-right minimum} of $\pi$ is a letter with no smaller
letter to the left of it; the number of left-to-right minima is
denoted $\lmin(\pi)$.  The statistics \emph{right-to-left minima}
($\rmin$), \emph{left-to-right maxima} ($\lmax$), and
\emph{right-to-left maxima} ($\rmax$) are defined similarly.

In the following table we define the remaining statistics that are of
interest to us. For reference we include the
statistics already defined in the preceding few paragraphs.

\begin{list}{}{
    \setlength{\leftmargin}{5em}
    \setlength{\itemsep}{0.5ex}
    \setlength{\labelwidth}{7em}
  }

\item[$\asc\, = $] number of ascent;

\item[$\comp\, = $] number of components;

\item[$\des\, = $] number of descents;

\item[$\exc\, = $] number of excedances: positions $i$ such that $a_i>i$;

\item[$\fp\, = $] number of fixed points: positions $i$ such that $a_i=i$;

\item[$\head\, = $] first element: $\head(\pi)=a_1$;

\item[$\last\, = $] last element: $\last(\pi)=a_n$;

\item[$\ldr\, = $] length of the leftmost decreasing run: largest $i$
  such that $a_1>a_2>\dots>a_i$;

\item[$\lds\, = $] length of the longest decreasing sequence in a
  permutation;

\item[$\lir\, = $] length of the leftmost increasing run: largest $i$
  such that $a_1<a_2<\dots<a_i$;

\item[$\lis\, = $] length of the longest increasing sequence in a
  permutation;

\item[$\lmax\, = $] number of left-to-right maxima; 

\item[$\lmin\, = $] number of left-to-right minima;

\item[$\peak\, = $] number of peaks: positions $i$ in $\pi$ such that
  $a_{i-1}<a_i>a_{i+1}$;

\item[$\rank\, = $] largest $k$ such that $a_i>k$ for all $i\leq k$ (see
  \cite{EP04});

\item[$\rdr\, = $] $\lirR\,=\,$ length of the rightmost decreasing run;

\item[$\rmax\, = $] number of right-to-left maxima;

\item[$\rmin\, = $] number of right-to-left minima;

\item[$\rir\, = $] $\ldrR \,=\,$ length of the rightmost increasing run;

\item[$\slmax\, = $] the number of letters to the left of second
  left-to-right maximum in $\pi\infty$: largest $i$ such that $a_1\geq
  a_1$, $a_1\geq a_2$, \dots, $a_1\geq a_i$;

\item[$\valley\, = $] number of valleys: positions $i$ in $\pi$ such that
  $a_{i-1}>a_i<a_{i+1}$;

\item[$\zeil\, = $] $\rdrI\,=\,$ length of the longest subword
  $n(n-1)\dots i$ (see \cite{Z92}).

\end{list}\smallskip

Let us also describe some of the derived statistics:

\begin{list}{}{
    \setlength{\leftmargin}{5em}
    \setlength{\itemsep}{0.5ex}
    \setlength{\labelwidth}{7em}
  }

\item[$\compR\, = $] number of reverse components: a \emph{reverse
  component} is a nonempty segment $\tau$ of $\pi$ such that $\pi =
  \sigma\tau\rho$ with $\sigma\succ\tau\succ\rho$, and such that if
  $\tau$ = $\alpha\beta$ and $\alpha\succ\beta$ then $\alpha$ or
  $\beta$ is empty;

\item[$\headI\, = $] position of the smallest letter;

\item[$\lastI\, = $] position of the largest letter;

\item[$\lirI\, = $] $\zeilC \,=\,$ largest $i$ such that $12\dots i$ is a
  subword in $\pi$;

\item[$\peakI\, = $] number of letters $a_i$ that are to the right of both
  $a_i-1$ and $a_i+1$;

\item[$\valleyI\, = $] number of letters $a_i$ that are to the left of
  both $a_i-1$ and $a_i+1$.

\end{list}

\section{Bijections in the literature}
\label{sec:bij}

In this section we describe the bijections, and we try to stay close
to the original sources when doing so. 
In what follows $\sym_n(\tau)$ is the set of $\tau$-avoiding
permutations of length $n$, and $\dyck_n$ is the set of Dyck paths of
length $2n$.

\subsection{Knuth's bijection, 1973}\label{subsec:Knuth}

Knuth~\cite[pp. 242--243]{K69} gives a bijection from $312$-avoiding
permutations to ``stack words''. Formulated a bit differently, it
amounts to a bijection from $132$-avoiding permutations to Dyck
paths. Knuth~\cite[pp. 60--61]{K73} also gives a bijection from
$321$-avoiding permutations to Dyck paths. By letting permutations that
are mapped to the same Dyck path correspond to each other, a bijection
between $321$- and $132$-avoiding permutation is obtained---we call it
Knuth's bijection.

We start by describing the bijection from $132$-avoiding permutations
to Dyck paths. We shall refer to it as the {\em standard bijection}.
(This bijection is the same as the one given by
Krattenthaler~\cite{K01}, who, however, gives a non-recursive
description of it; see Section~\ref{subsec:Krattan}.) Let $\pi=\pi_Ln\pi_R$
be a $132$-avoiding permutation of length $n$. Each letter of $\pi_L$
is larger than any letter of $\pi_R$, or else a $132$ pattern would
be formed. We define the standard bijection $f$ recursively by
$f(\pi)=uf(\pi_L)df(\pi_R)$ and $f(\emptyword)=\emptyword$. Here, and
elsewhere, $\emptyword$ denotes the empty word/permutation. Thus,
under the standard bijection, the position of the largest letter in a
$132$-avoiding permutation determines the first return to $x$-axis and
vice versa. For instance,
\begin{align*}
f(7564213) 
&= udf(564213) 
= uduf(5)df(4213) 
= uduuddudf(213) \\
&= uduudduduf(21)d
= uduudduduudf(1)d \\
&= uduudduduududd \\
&= 
\vcenter{\xymatrix@R=.4ex@C=.4ex@!{
      &      &      &      &\n\erd &      &      &      &       &      &\n\erd &      &\n\erd &      &   \\
      &\n\erd &      &\n\eru &      &\n\erd &      &\n\erd &       &\n\eru &      &\n\eru &      &\n\erd &   \\
\n\eru &      &\n\eru &      &      &      &\n\eru &      & \n\eru &      &      &      &      &      &\n 
}}
\end{align*}

As mentioned, Knuth also gives a bijection from $321$-avoiding
permutations to Dyck paths: Given a $321$-avoiding permutation, start
by applying the {\em Robinson-Schensted-Knuth correspondence} to
it. This classic correspondence gives a bijection between permutations
$\pi$ of length $n$ and pairs $(P,Q)$ of {\em standard Young tableaux}
of the same shape $\lambda\vdash n$. As is well known, the length of
the longest decreasing subword in $\pi$ corresponds to the number
of rows in $P$ (or $Q$). Thus, for $321$-avoiding permutations, the
tableaux $P$ and $Q$ have at most two rows.

The {\em insertion tableau} $P$ is obtained by reading
$\pi=a_1a_2\dots a_n$ from left to right and, at each step, inserting
$a_i$ to the partial tableau obtained thus far. Assume that
$a_1,a_2,\dots,a_{i-1}$ have been inserted. If $a_i$ is larger than
all the elements in the first row of the current tableau, place $a_i$
at the end of the first row. Otherwise, let $m$ be the leftmost
element in the first row that is larger than $a_i$. Place $a_i$ in the
square that is occupied by $m$, and place $m$ at the end of the second
row.  The {\em recording tableau} $Q$ has the same shape as $P$ and is
obtained by placing $i$, for $i$ from $1$ to $n$, in the position of
the square that in the construction of $P$ was created at step $i$
(when $a_i$ was inserted). For example, the pair of tableaux
corresponding to the $321$-avoiding permutation $3156247$ we get by the
following sequence of insertions:
\newcommand{\f}[1]{\makebox[3.4ex][l]{#1}}
\begin{align*}
\left(\!
\begin{array}{l|l}
  \emptyword & \emptyword \\
\end{array}
\!\right)
\to
\left(\!
\begin{array}{l|l}
  3 & 1 \\
\end{array}
\!\right)
&\to
\left(\!
\begin{array}{l|l}
  \f 1 & \f 1 \\
  \f 3 & \f 2
\end{array}
\!\right) \\
&
\to
\left(\!
\begin{array}{l|l}
  \f{15} & \f{13} \\
  \f 3   & \f 2
\end{array}
\!\right)
\to
\left(\!
\begin{array}{l|l}
  \f{156} & \f{134} \\
  \f 3    & \f 2
\end{array}
\!\right) \\
& \to
\left(\!
\begin{array}{l|l}
  \f{126} & \f{134} \\
  \f{35}  & \f{25}
\end{array}
\!\right)
\to
\left(\!
\begin{array}{l|l}
  \f{124} & \f{134} \\
  \f{356} & \f{256}
\end{array}
\!\right)
\to
\left(\!
\begin{array}{l|l}
  1247 & 1347 \\
  356  & 256
\end{array}
\!\right).
\end{align*}

The pair of tableaux $(P,Q)$ is then turned into a Dyck path $D$. The
first half, $A$, of the Dyck path we get by recording, for $i$ from
$1$ to $n$, an up-step if $i$ is in the first row of $P$, and a
down-step if it is in the second row. Let $B$ be the word obtained
from $Q$ in the same way but interchanging the roles of $u$ and
$d$. Then $D=AB^r$ where $B^r$ is the reverse of $B$. Continuing with
the example above we get
$$D =
\vcenter
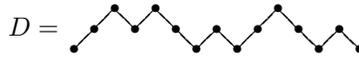
{\xymatrix@R=.4ex@C=.4ex@!{
      &      &\n\erd &      &\n\erd &      &      &      &      &      &\n\erd & \\
      &\n\eru &      &\n\eru &      &\n\erd &      &\n\erd &      &\n\eru &      &\n\erd &      &\n\erd     \\
\n\eru &      &      &      &      &      &\n\eru &      &\n\eru &      &      &      &\n\eru &      &\n
}}
$$

Elizalde and Pak~\cite{EP04} use this bijection together with a
slight modification of the standard bijection to give a
combinatorial proof of a generalization of the result by Robertson et
al.~\cite{RSZ03} that fixed points have the same distribution on
$123$- and $132$-avoiding permutations. The modification they use is
to reflect the Dyck path obtained from the standard bijection with
respect to the vertical line crossing the path in the middle.
Alternatively, the path can be read from the permutation diagram as
described in~\cite{EP04}. We follow Elizalde and Pak and apply the
same modification. After reflection, the path $f(7564213)$ above is
the same as the path $D$ in the preceding example. Thus the image of
the $321$-avoiding permutation $3156247$ under what we call Knuth's
bijection is the $132$-avoiding permutation $7564213$.

\subsection{Knuth-Rotem's bijection, 1975}

Rotem~\cite{R75} gives a bijection between 321-avoiding permutations
and Dyck paths, described below. Combining it with the standard
bijection gives a bijection from $321$- to $132$-avoiding
permutations---we call it {\em Knuth-Rotem's bijection}.

A {\em ballot-sequence} $b_1b_2\dots b_n$ satisfies the two
conditions
\begin{enumerate}
\item $b_1\leq b_2\leq \dots\leq b_n$;
\item $0\leq b_i\leq i-1$, for $i=1,2,\dots,n$.
\end{enumerate}

Let $\pi=p_1 p_2\dots p_n$ be a $321$-avoiding permutation. From it
we construct a ballot-sequences $b_1b_2\dots b_n$: Let $b_1=0$. For
$i=2,\dots,n$, let $b_i=b_{i-1}$ if $ p_i$ is a left-to-right maximum
in $\pi$, and let $b_i= p_i$ otherwise. 

For the permutation $\pi=2513476$ we get the ballot-sequences
$0013446$. This sequence we represent by a ``bar-diagram'', which
in turn can be viewed as a lattice path from $(0,0)$ to $(7,7)$:
$$
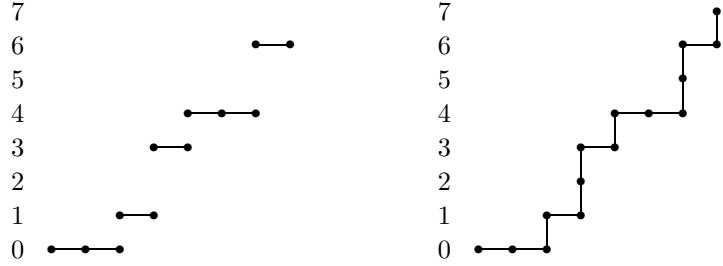
\begin{array}{cc}
\xymatrix@R=.1ex@C=.1ex@!{
7 \\
6 &       &       &       &      &      &      &\n\er &\n \\
5 &       &       &       &      &      &      &      & \\
4 &       &       &       &      &\n\er &\n\er &\n \\
3 &       &       &       &\n\er &\n \\
2 &       &       &       &      &  \\
1 &       &       & \n\er & \n\\
0 & \n\er & \n\er & \n \\
}\qquad
&\qquad
\xymatrix@R=.1ex@C=.1ex@!{
7 &       &       &       &      &      &      &      &\n \\
6 &       &       &       &      &      &      &\n\er &\n\eu \\
5 &       &       &       &      &      &      &\n\eu & \\
4 &       &       &       &      &\n\er &\n\er &\n\eu \\
3 &       &       &       &\n\er &\n\eu \\
2 &       &       &       &\n\eu &  \\
1 &       &       & \n\er &\n\eu\\
0 & \n\er & \n\er & \n\eu \\
}
\end{array}
$$
Rotating that path counter clockwise by $3\pi/4$ radians we get  
$$
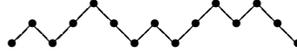
\vcenter{\xymatrix@R=.4ex@C=.4ex@!{
      &      &      &      &\n\erd &      &      &      &       &      &\n\erd &      &\n\erd &      &   \\
      &\n\erd &      &\n\eru &      &\n\erd &      &\n\erd &       &\n\eru &      &\n\eru &      &\n\erd &   \\
\n\eru &      &\n\eru &      &      &      &\n\eru &      & \n\eru &      &      &      &      &      &\n 
}}
$$ In the previous subsection we saw that this path is $f(7564213)$
where $f$ is the standard bijection. Thus the image of the
$321$-avoiding permutation $2513476$ under Knuth-Rotem's bijection is
the $132$-avoiding permutation $7564213$.

\subsection{Simion-Schmidt's bijection, 1985} 
Consider the following algorithm:\medskip

\begin{tabular}{ll}
Input:  & A permutation $\sigma=a_1a_2\dots,a_n$ in $\sym_n(123)$. \\
Output: & A permutation  $\tau=c_1c_2\dots c_n$ in $\sym_n(132)$. \bigskip\\
1 &\mbox{} $c_1:=a_1;$ $x:=a_1$\medskip\\
2 &\mbox{} {\tt for} $i=2, \dots,  n$: \\
3 &\mbox{} \quad\; {\tt if} $a_i < x$: \\
4 &\mbox{} \quad\;\quad\; $c_i:=a_i$; $x:=a_i$\\
5 &\mbox{} \quad\; {\tt else}:\\
6 &\mbox{} \quad\;\quad\;
$c_i:=\min\{\,k\ |\ x<k\leq n,\,k\neq c_j\textrm{ for all } j<i\,\}$\\
\end{tabular}\bigskip

\noindent
The map $\sigma\mapsto\tau$ is the Simion-Schmidt
bijection~\cite{SS85}. As an example, the $123$-avoiding permutation
$6743152$ maps to the $132$-avoiding permutation $6743125$.

\subsection{Knuth-Richards' bijection, 1988}\label{subsec:KRi}
Richards' bijection~\cite{R88} from Dyck paths to $123$-avoiding
permutations is given by the following algorithm:\medskip

\begin{tabular}{ll}
Input:  & A Dyck path $P=b_1b_2\dots b_{2n}$. \\
Output: & A permutation $\pi=a_1 a_2\dots a_n$ in $\sym_n(123)$.\bigskip\\
1 &\mbox{} $r:=n+1$; $s:=n+1$; $j:=1$ \medskip\\
2 &\mbox{} {\tt for} $i=1, \dots, n$:\\
3 &\mbox{} \quad\; {\tt if} $b_j$ is an up-step: \\
4 &\mbox{} \quad\;\quad\; 
{\tt repeat} $s:=s-1$; $j:=j+1$ {\tt until} $b_j$ is a down-step \\  
5 &\mbox{} \quad\;\quad\; $a_s:=i$ \\
6 &\mbox{} \quad\; {\tt else:} \\
7 &\mbox{} \quad\;\quad\; 
{\tt repeat} $r:=r-1$ {\tt until} $a_r$ is unset\\
8 &\mbox{} \quad\;\quad\; $a_r:=i$ \\
9 &\mbox{} \quad\; $j:=j+1$ 
\end{tabular}\bigskip

\noindent
The Knuth-Richards bijection, from $\sym_n(132)$ to $\sym_n(123)$, is
defined by
$$\KRi \ =\ \Ri\,\circ\, f,
$$ 
where $f$ is the standard bijection from $132$-avoiding
permutations to Dyck paths, and $\Ri$ is the algorithm just
described. As an example, applying Knuth-Richards' bijection to
$6743125$ yields $5743612$.

\subsection{West's bijection, 1995}\label{subsec:west}

West's bijection~\cite{W95} is induced by an isomorphism between {\em
  generating trees}. The two isomorphic trees generate 123- and
132-avoiding permutations, respectively. We give a brief description
of that bijection: Given a permutation $\pi=p_1 p_2\dots p_{n-1}$ and
a positive integer $i\leq n$, let
$$\pi^i\ =\ p_1\, \dots\,  p_{i-1}\ n\ p_i\,\dots\, p_{n-1};
$$ we call this {\em inserting $n$ into site $i$}. With respect to a
fixed pattern $\tau$ we call site $i$ of $\pi$ in $\sym_{n-1}(\tau)$
{\em active} if the insertion of $n$ into site $i$ creates a
permutation in $\sym_n(\tau)$.

For $i=0,\dots, n-1$, let $a_{i+1}$ be the number of active sites in
the permutation obtained from $\pi$ by removing the $i$ largest
letters. The {\em signature} of $\pi$ is the word 
$$a_0a_1\dots a_{n-1}.
$$ 
West~\cite{W95} showed that for $123$-avoiding permutations, as
well as for $132$-avoiding permutations, the signature determines the
permutation uniquely. This induces a natural bijection between the two
sets. For example, the $123$-avoiding permutation $536142$ corresponds
to the $132$-avoiding permutation $534612$---both have
the same signature, $343322$.

\subsection{Krattenthaler's bijection, 2001}\label{subsec:Krattan}

Krattenthaler's bijection~\cite{K01} uses Dyck paths as intermediate
objects. Permutations that are mapped to the same Dyck path correspond
to each other under this bijection.

The first part of Krattenthaler's bijection is a bijection from
$123$-avoiding permutations to Dyck paths. Reading right to left,
let the right-to-left maxima in $\pi$ be $m_1$, $m_2$, \dots, $m_s$, so
that 
$$\pi \,=\, w_sm_s\dots w_2m_2 w_1m_1,
$$ where $w_i$ is the subword of $\pi$ in between $m_{i+1}$ and
$m_i$. Since $\pi$ is $123$-avoiding, the letters in $w_i$ are in
decreasing order.  Moreover, all letters of $w_i$ are smaller than
those of $w_{i+1}$.

To define the bijection, read $\pi$ from right to left. Any
right-to-left maximum $m_i$ is translated into $m_i-m_{i-1}$ up-steps
(with the convention $m_0=0$). Any subword $w_i$ is translated into
$|w_i|+1$ down-steps, where $|w_i|$ denotes the number of letters of
$w_i$. Finally, the resulting path is reflected in a vertical line
through the center of the path. Alternatively, we could have
generated the Dyck path from right to left.

The second part of Krattenthaler's bijection is a bijection from
$132$-avoiding permutations to Dyck paths. Read $\pi= p_1 p_2\dots
p_n$ in $\sym_n(132)$ from left to right and generate a Dyck
path. When $ p_j$ is read, adjoin, to the path obtained thus far, as
many up-steps as necessary to reach height $h_j+1$, followed by a
down-step to height $h_j$ (measured from the $x$-axis); here $h_j$ is
the number of letters in $ p_{j+1}\dots p_n$ which are larger than $
p_j$. This procedure can be shown to be equivalent to the standard
bijection from $132$-avoiding permutations to Dyck paths.

For instance, Krattenthaler's bijection sends the permutation $536142$
in $\sym_6(123)$ to the permutation $452316$ in $\sym_6(132)$---both
map to the same Dyck path,
$$
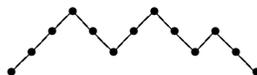
\vcenter{\xymatrix@R=.4ex@C=.4ex@!{
       &       &       &\n\erd &       &       &       &\n\erd &       &       &       &       & \\
       &       &\n\eru &       &\n\erd &       &\n\eru &       &\n\erd &       &\n\erd &       & \\
       &\n\eru &       &       &       &\n\eru &       &       &       &\n\eru &       &\n\erd & \\
\n\eru &       &       &       &       &       &       &       &       &       &       &       &\n 
}}
$$


\subsection{Reifegerste's bijection, 2002}\label{subsec:astrid}
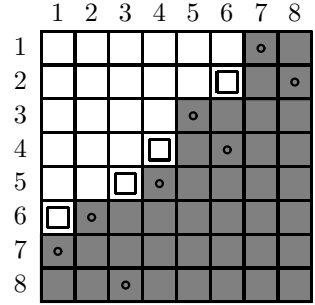
$$
\setlength{\unitlength}{0.9mm}
\begin{picture}(48,44)(-15,1)
\thicklines
\put(0,0){ \put(-10,0){\shade\path(0,0)(40,0)(40,10)(0,10)(0,0)}
\put(-10,0){\shade\path(5,10)(40,10)(40,15)(5,15)(5,10)}
\put(-10,0){\shade\path(15,15)(40,15)(40,20)(15,20)(15,15)}
\put(-10,0){\shade\path(20,20)(40,20)(40,30)(20,30)(20,20)}
\put(-10,0){\shade\path(30,30)(40,30)(40,40)(30,40)(30,30)}
\put(-10,0){\grid(40,40)(5,5)} \put(-9,11){\grid(3,3)(3,3)}
\put(1,16){\grid(3,3)(3,3)}
\put(6,21){\grid(3,3)(3,3)}\put(16,31){\grid(3,3)(3,3)}
\put(22.5,37.5){\circle{1}}\put(27.5,32.5){\circle{1}}\put(12.5,27.5){\circle{1}}
\put(17.5,22.5){\circle{1}}\put(7.5,17.5){\circle{1}}
\put(-2.5,12.5){\circle{1}}\put(-7.5,7.5){\circle{1}}\put(2.5,2.5){\circle{1}}
\put(-14,1.5){8}\put(-14,6.5){7}\put(-14,11.5){6}\put(-14,16.5){5}
\put(-14,21.5){4}\put(-14,26.5){3}\put(-14,31.5){2}\put(-14,36.5){1}
\put(-8.5,41.5){1}\put(-3.5,41.5){2}\put(1.5,41.5){3}\put(6.5,41.5){4}
\put(11.5,41.5){5}\put(16.5,41.5){6}\put(21.5,41.5){7}\put(26.5,41.5){8}}
\end{picture}\;\;
\parbox{54ex}{\vspace{-23ex} This figure illustrates Reifegerste's
  bijection~\cite{R02}. It pictures the $321$-avoiding permutation
  $\pi=13256847$ and the $132$-avoiding permutation $\pi'=78564213$,
  two permutations that correspond to each other under that bijection.

  Let $\pi= a_1 a_2\dots a_n$ be a $321$-avoiding permutation,
  and let $E$ be the set of pairs\medskip \\ 
  $\mbox{}\qquad\quad E\,=\,\{\,(i, a_i)\ |\ i \textrm{ is an excedance}\,\}.
  $ 
  }
$$ For each pair $(i,a_i)$ in $E$, we place a square, called an {\em
  $E$-square}, in position $(i,n+1-a_i)$ in an $n\times n$ {\em
  permutation matrix}. ($E$ uniquely determines $\pi$.) Next we shade
each square $(a,b)$ of the matrix where there are no $E$-squares in
the region $\{(i,j)\ |\ i\geq a,\ j\geq b\}$, thus obtaining a {\em
  Ferrer's diagram}. Finally, we get the $132$-avoiding permutation
$\pi'$ corresponding to $\pi$ by placing dots (circles), row by row
starting from the first row, in the leftmost available shaded square
such that there are no two dots in any column or row. If $(i,j)$
contains a dot, then $\pi'(i)=j$.

\subsection{Elizalde-Deutsch's bijection, 2003}

Here is an outline of a bijection by Elizalde and Deutsch~\cite{ED03}:
Map $321$- and $132$-avoiding permutation bijectively to Dyck paths;
use an automorphism $\Psi$ on Dyck paths; and match permutations with
equal paths.

We start by describing the automorphism $\Psi$. Let $P$ be a Dyck path
of length $2n$. Each up-step of $P$ has a corresponding down-step in the
sense that the path between the up-step and the down-step form a
proper Dyck path. Match such pairs of steps. Let $\sigma$ in $\sym_{2n}$
be the permutation defined by $\sigma_i=(i+1)/2$ if $i$ is odd, and
$\sigma_i=2n+1-i/2$ otherwise. For $i$ from $1$ to $2n$, consider the
$\sigma_i$-th step of $P$. If the corresponding matching step has not
yet been read, define the $i$-th step of $\Psi(P)$ to be an up-step,
otherwise let it be a down-step. For example,
$$\Psi(uuduudududddud) \,=\, uuuddduduuddud.
$$

The bijection $\psi$ from $321$-avoiding permutations to $\dyck_n$ is
defined as follows. Any permutation $\pi$ in $\sym_n$ can be
represented as an $n\times n$ array with crosses in the squares
$(i,\pi(i))$. Given the array of $\pi$ in $\sym_n(321)$, consider the
path with {\em down-} and {\em right-}steps along the edges of the
squares that goes from the upper-left corner to the lower-right
corner of the array leaving all the crosses to the right and
remaining always as close to the main diagonal as possible. Then the
corresponding Dyck path is obtained from this path by reading an
up-step every time the path moves down, and a down-step every time
the path moves to the right. For example, 
$$\psi(2314657) \,=\, uuduudududddud.
$$

The bijection $\phi$ from $132$-avoiding permutations to $\dyck_n$ is the
standard bijection followed by a reflection of the path with respect
to a vertical line through the middle of the path. For example,
$$\phi(7432516) \,=\, uuduudududddud.
$$

The Elizalde-Deutsch bijection, from $\sym_n(321)$ to $\sym_n(132)$,
is defined by
$$\ED \ =\ \phi^{-1}\circ\,\Psi^{-1}\circ\,\psi.
$$
As an example, it send $2314657$ to $2314657$.

\subsection{Mansour-Deng-Du's bijection, 2006}
Let $i$ be a positive integer smaller than $n$.  Let
$s_i:\sym_n\rightarrow\sym_n$ act on permutations by interchanging the
letters in positions $i$ and $i+1$.  We call $s_i$ a {\em simple
  transposition}, and write the action of $s_i$ as $\pi s_i$. So,
$\pi(s_is_j)=(\pi s_i)s_j$. For any permutation $\pi$ of length $n$,
the {\em canonical reduced decomposition} of $\pi$ is
$$\pi=(12\dots n)\sigma=(12\dots n)\sigma_1\sigma_2\dots\sigma_k,
$$ where $\sigma_i=s_{h_i}s_{h_i-1}\dots s_{t_i}$, $h_i\geq t_i$, 
$1\leq i\leq k$ and $1\leq h_1<h_2<\dots <h_k\leq n-1$. For example,
$415263=(s_3s_2s_1)(s_4s_3)(s_5)$. 

Mansour, Deng and Du~\cite{MDD} use canonical reduced decompositions
to construct a bijection between $\sym_n(321)$ and $\sym_n(231)$. They
show that a permutation is $321$-avoiding precisely when $t_i\geq
t_{i-1}+1$ for $2\leq i\leq k$~\cite[Thm. 2]{MDD}. They also show that
a permutation is $231$-avoiding precisely when $t_i\geq t_{i-1}$ or
$t_i\geq h_{i-j}+2$ for $2\leq i\leq k$ and $1\leq j\leq
i-1$~\cite[Thm. 15]{MDD}. Using these two theorems they build their
bijection, which is composed of two bijections: one from $\sym_n(321)$
to $\dyck_n$, and one from $\sym_n(231)$ to $\dyck_n$.

For a Dyck path $P$, we define the $(x+y)$-labeling of $P$ as
follows: each cell in the region enclosed by $P$ and the $x$-axis,
whose corner points are $(i,j)$, $(i+1,j-1)$, $(i+2,j)$ and
$(i+1,j+1)$ is labeled by $(i+j)/2$. If $(i-1,j-1)$ and $(i,j)$ are
starting points of two consecutive up-steps, then we call the cell
with leftmost corner $(i,j)$ an {\em essential cell} and the
up-step $((i-1,j-1),(i,j))$ its {\em left arm}. We define the {\em
  zigzag strip} of $P$ as follows: If there is no essential cell in
$P$, then the zigzag strip is simply the empty set. Otherwise, we
define the zigzag strip of $P$ as the border strip that begins at the
rightmost essential cell.  For example, the zigzag strip of the Dyck
path $uuduuududddudduduuddud$ in Figure~\ref{dyck01} is the shaded
cell labeled by 9, while for the Dyck path $uuduuududdd$ (obtained
from that in Figure~\ref{dyck01} by ignoring the steps $15$ to $22$)
the zigzag strip is the shaded connected cells labeled by $2$, $3$,
$4$, $5$ and $6$.

\begin{figure}
\begin{center}

\setlength{\unitlength}{5mm}

\begin{picture}(23,5)
\thicklines
\shade\path(5,3)(6,4)(7,3)(8,4)(9,3)(10,2)(11,1)(12,2)(13,1)(12,0)(11,1)(10,0)(9,1)(8,2)(7,3)(6,2)(5,3)

\shade\path(1,1)(2,2)(3,1)(4,2)(5,1)(6,2)(7,1)(6,0)(5,1)(4,0)(3,1)(2,0)(1,1)

\path(0,0)(1,1)(2,2)(3,1)(4,2)(5,3)(6,4)(7,3)(8,4)(9,3)(10,2)(11,1)(12,2)(13,1)(14,0)(15,1)(16,0)(17,1)
(18,2)(19,1)(20,0)(21,1)(22,0)

\path(1,1)(2,0)(3,1)(4,0)(5,1)(6,0)(7,1)(8,0)(9,1)(10,0)(11,1)(12,0)(13,1)(14,0)(15,1)(16,0)(17,1)(18,0)(19,1)(20,0)

\path(4,2)(5,1)(6,2)(7,1)(8,2)(9,1)(10,2)

\path(5,3)(6,2)(7,3)(8,2)(9,3)

\shade\path(17,1)(18,2)(19,1)(18,0)

\put(1.8,0.7){\large 1}\put(3.8,0.7){\large 2}\put(5.8,0.7){\large
3}\put(7.8,0.7){\large 4}\put(9.8,0.7){\large
5}\put(11.8,0.7){\large 6} \put(17.8,0.7){\large 9}

\put(4.8,1.7){\large 3}\put(6.8,1.7){\large 4}\put(8.8,1.7){\large
5}

\put(5.8,2.7){\large 4}\put(7.8,2.7){\large 5}

\path(0,5)(0,0)(23,0)

\path(-0.2,4.7)(0,5)(0.2,4.7)
\path(22.7,0.2)(23,0)(22.7,-0.2)

\path(1,-0.1)(1,0.1)\path(3,-0.1)(3,0.1)\path(5,-0.1)(5,0.1)\path(7,-0.1)(7,0.1)\path(9,-0.1)(9,0.1)
\path(11,-0.1)(11,0.1)\path(13,-0.1)(13,0.1)\path(15,-0.1)(15,0.1)\path(17,-0.1)(17,0.1)\path(19,-0.1)(19,0.1)
\path(21,-0.1)(21,0.1)
\path(-0.1,1)(0.1,1)\path(-0.1,2)(0.1,2)\path(-0.1,3)(0.1,3)\path(-0.1,4)(0.1,4)

\put(-0.6,0.8){\small 1} \put(-0.6,1.8){\small
2}\put(-0.6,2.8){\small 3} \put(-0.6,3.8){\small 4}

\put(-0.15,-0.7){\small 0}\put(0.85,-0.7){\small
1}\put(1.85,-0.7){\small 2} \put(2.85,-0.7){\small
3}\put(3.85,-0.7){\small 4}\put(4.85,-0.7){\small
5}\put(5.85,-0.7){\small 6}\put(6.85,-0.7){\small
7}\put(7.85,-0.7){\small 8}\put(8.85,-0.7){\small
9}\put(9.75,-0.7){\small 10}\put(10.75,-0.7){\small
11}\put(11.75,-0.7){\small 12}\put(12.75,-0.7){\small
13}\put(13.75,-0.7){\small 14}\put(14.75,-0.7){\small
15}\put(15.75,-0.7){\small 16}\put(16.75,-0.7){\small
17}\put(17.75,-0.7){\small 18}\put(18.75,-0.7){\small
19}\put(19.75,-0.7){\small 20}\put(20.75,-0.7){\small
21}\put(21.75,-0.7){\small 22}

\end{picture}
\caption{}\label{dyck01}
\smallskip
\end{center}
\end{figure}

Let $P_{n,k}$ be a Dyck path of semi-length $n$ containing $k$
essential cells. We define its {\em zigzag decomposition} as follows:
The zigzag decomposition of $P_{n,0}$ is the empty set.  The zigzag
decomposition of $P_{n,1}$ is the zigzag strip.  If $k\geq 2$, then we
decompose $P_{n,k}=P_{n,k-1}Q$, where $Q$ is the zigzag strip of
$P_{n,k}$ and $P_{n,k-1}$ is the Dyck path obtained from $P$ by
deleting $Q$. Reading the labels of $Q$ from left to right, ignoring
repetitions, we get a sequence of numbers $\{i,i+1,\dots,j\}$, and we
  associate $Q$ with the sequence of simple decompositions
  $\sigma_k=s_js_{j-1}\dots s_i$.  For $P_{n,i}$ with $i\leq k-1$ repeat
  the above procedure to get $\sigma_{k-1}$, $\dots$, $\sigma_2$,
  $\sigma_1$.  The zigzag decomposition of $P_{n,k}$ is then given by
  $\sigma=\sigma_1\sigma_2\dots\sigma_k$.

From the zigzag decomposition we get a $321$-avoiding permutation
$\pi=(12\dots n)\sigma$ whose canonical reduced decomposition is
$\sigma$. For the Dyck path $P_{11,4}$
in Figure~\ref{dyck01} we have
$$\sigma=(s_3s_2s_1)(s_4s_3)(s_6s_5s_4)(s_9)
$$ 
and the corresponding
permutation in $\sym_{11}(321)$ is $(4,1,5,7,2,3,6,8,10,9,11)$. 

\begin{figure}
\begin{center}

\setlength{\unitlength}{5mm}

\begin{picture}(23,5)
\thicklines

\shade\path(5,3)(6,4)(7,3)(8,4)(9,3)(8,2)(7,3)(6,2)(5,3)

\shade\path(1,1)(2,2)(3,1)(4,2)(5,1)(6,2)(7,1)(8,2)(9,1)(10,2)(11,1)(12,2)(13,1)
(12,0)(11,1)(10,0)(9,1)(8,0)(7,1)(6,0)(5,1)(4,0)(3,1)(2,0)(1,1)

\path(0,0)(1,1)(2,2)(3,1)(4,2)(5,3)(6,4)(7,3)(8,4)(9,3)(10,2)(11,1)(12,2)(13,1)(14,0)(15,1)(16,0)(17,1)
(18,2)(19,1)(20,0)(21,1)(22,0)

\path(1,1)(2,0)(3,1)(4,0)(5,1)(6,0)(7,1)(8,0)(9,1)(10,0)(11,1)(12,0)(13,1)(14,0)(15,1)(16,0)(17,1)(18,0)(19,1)(20,0)

\path(4,2)(5,1)(6,2)(7,1)(8,2)(9,1)(10,2)

\path(5,3)(6,2)(7,3)(8,2)(9,3)

\shade\path(17,1)(18,2)(19,1)(18,0)

\put(1.8,0.7){\large 1}\put(3.8,0.7){\large 2}\put(5.8,0.7){\large
3}\put(7.8,0.7){\large 4}\put(9.8,0.7){\large
5}\put(11.8,0.7){\large 6} \put(17.8,0.7){\large 9}

\put(4.8,1.7){\large 2}\put(6.8,1.7){\large 3}\put(8.8,1.7){\large
4}

\put(5.8,2.7){\large 2}\put(7.8,2.7){\large 3}

\path(0,5)(0,0)(23,0)
\path(-0.2,4.7)(0,5)(0.2,4.7)
\path(22.7,0.2)(23,0)(22.7,-0.2)

\path(1,-0.1)(1,0.1)\path(3,-0.1)(3,0.1)\path(5,-0.1)(5,0.1)\path(7,-0.1)(7,0.1)\path(9,-0.1)(9,0.1)
\path(11,-0.1)(11,0.1)\path(13,-0.1)(13,0.1)\path(15,-0.1)(15,0.1)\path(17,-0.1)(17,0.1)\path(19,-0.1)(19,0.1)
\path(21,-0.1)(21,0.1)
\path(-0.1,1)(0.1,1)\path(-0.1,2)(0.1,2)\path(-0.1,3)(0.1,3)\path(-0.1,4)(0.1,4)

\put(-0.6,0.8){\small 1} \put(-0.6,1.8){\small
2}\put(-0.6,2.8){\small 3} \put(-0.6,3.8){\small 4}

\put(-0.15,-0.7){\small 0}\put(0.85,-0.7){\small
1}\put(1.85,-0.7){\small 2} \put(2.85,-0.7){\small
3}\put(3.85,-0.7){\small 4}\put(4.85,-0.7){\small
5}\put(5.85,-0.7){\small 6}\put(6.85,-0.7){\small
7}\put(7.85,-0.7){\small 8}\put(8.85,-0.7){\small
9}\put(9.75,-0.7){\small 10}\put(10.75,-0.7){\small
11}\put(11.75,-0.7){\small 12}\put(12.75,-0.7){\small
13}\put(13.75,-0.7){\small 14}\put(14.75,-0.7){\small
15}\put(15.75,-0.7){\small 16}\put(16.75,-0.7){\small
17}\put(17.75,-0.7){\small 18}\put(18.75,-0.7){\small
19}\put(19.75,-0.7){\small 20}\put(20.75,-0.7){\small
21}\put(21.75,-0.7){\small 22}

\end{picture}
\caption{}\label{dyck02}
\end{center}
\end{figure}

We will now describe a map from Dyck paths to 231-avoiding
permutations. For a Dyck path $P$, we define $(x-y)$-labeling of $P$
as follows (this labeling seems to be considered for the first time
in~\cite{BK}): each cell in the region enclosed by $P$ and the
$x$-axis, whose corner points are $(i,j)$, $(i+1,j-1)$, $(i+2,j)$ and
$(i+1,j+1)$ is labeled by $(i-j+2)/2$. We define the {\em trapezoidal
  strip} of $P$ as follows: If there is no essential cell in $P$, then
the trapezoidal strip is simply the empty set. Otherwise, we define
the trapezoidal strip of $P$ as the horizontal strip that touches the
$x$-axis and starts at the rightmost essential cell. For example, the
trapezoidal strip of the Dyck path $uuduuududddudduduuddud$ in
Figure~\ref{dyck02} is the shaded cell labeled by $9$, while for the
Dyck path $uuduuududdd$ (obtained from that in Figure~\ref{dyck02} by
ignoring the steps $15$ to $22$) the zigzag strip is the down-most
shaded strip with labels $1$, $2$, $3$, $4$, $5$ and $6$.

Let $P_{n,k}$ be a Dyck path of semi-length $n$ containing $k$ essential
cells. We define its {\em trapezoidal decomposition} as follows:
The trapezoidal decomposition of $P_{n,0}$ is the empty set.  The
trapezoidal decomposition of $P_{n,1}$ is the trapezoidal strip.  If
$k\geq 2$, then we decompose $P_{n,k}$ into $P_{n,k}=Q_1uQ_2d$, where
$u$ is the left arm of the rightmost essential cell that touches the
$x$-axis, $d$ is the last down step of $P_{n,k}$, and $Q_1$ and $Q_2$
carry the labels in $P_{n,k}$. Reading the labels of the
trapezoidal strip of $P_{n,k}$ from left to right we get a
sequence $\{i,i+1,\dots,j\}$, and we set $\sigma_k=s_js_{j-1}\dots
s_i$.  Repeat the above procedure for $Q_1$ and $Q_2$. Suppose the
trapezoidal decomposition of $Q_1$ and $Q_2$ are $\sigma^{\prime}$ and
$\sigma^{\prime\prime}$ respectively, then the trapezoidal
decomposition for $P_{n,k}$ is
$\sigma=\sigma^{\prime}\sigma^{\prime\prime}\sigma_k$.

From the trapezoidal decomposition we get a $231$-avoiding permutation
$\pi=(12\dots n)\sigma$ whose canonical reduced
decomposition is $\sigma$. For the Dyck path $P_{11,4}$
in Figure~\ref{dyck02} we have
$$\sigma=(s_3s_2)(s_4s_3s_2)(s_6s_5s_4s_3s_2s_1)(s_9)
$$ and the corresponding permutation in $\sym_{11}(231)$ is
$(7,1,5,4,2,3,6,8,10,9,11)$.

The two maps involving Dyck paths described in this subsection
induce a bijection from 321-avoiding to 231-avoiding permutations.

\section{A recursive description of the Knuth-Richards bijection}
\label{sec:CK}

We call a permutation $\pi$ {\em indecomposable} if $\comp(\pi)=1$;
otherwise we call $\pi$ {\em decomposable}. Equivalently, if we define
the sum $\oplus$ on permutations by $\sigma\oplus\tau = \sigma\tau'$,
where $\tau'$ is obtained from $\tau$ by adding $|\sigma|$ to each of
its letters, then a permutation is indecomposable if it cannot be
written as the sum of two nonempty permutations.

We shall describe, separately for $231$- and $321$-avoiding
permutations, how to generate the indecomposable permutations, thus
inducing a bijection we call \CK.

For a permutation of length $n$ to be $231$-avoiding everything to the
left of $n$ has to be smaller than anything to the right of
$n$. Clearly, if there is at least one letter to the left of $n$,
then the permutation is decomposable (everything to the right of $n$,
including $n$, would form the last component). Thus a $231$-avoiding
permutation of length $n$ is indecomposable if and only if it starts
with $n$.

To build an indecomposable $231$-avoiding permutation of length $n$
from a $231$-avoiding permutation of length $n-1$ we simply prepend
$n$. Let us call this map $\alpha$. For instance,
$\alpha(2134)=52134$.

Given $k$ indecomposable $231$-avoiding permutations $\pi_1$, $\pi_2$,
\dots, $\pi_k$, we build the corresponding permutation by summing:
$\pi_1\oplus \pi_2\oplus\dots\oplus \pi_k$. Given $k$ indecomposable
$321$-avoiding permutations $\pi_1, \pi_2,\dots,\pi_k$ we build the
corresponding permutation by summing in reverse order: $\pi_k\oplus
\pi_{k-1}\oplus\dots\oplus \pi_1$.

Here is how we build an indecomposable $321$-avoiding
permutation $\pi'$ of length $n$ from a $321$-avoiding permutation
$\pi$ of length $n-1$:
$$
\vcenter{\xymatrix@R=0.5ex@C=0ex@!{
\pi &=&&2 &4                  &1 &3 &5 &*=<3.5ex>[F]{7} &6 &*=<3.5ex>[F]{9}  &8 \\ 
    & & 2 &4 &*=<3.5ex>[F]{10}&1 &3 &5 &*=<3.5ex>[F]{7} &6 &*=<3.5ex>[F]{9}  &8 \\ 
\pi'&=& 2 &4 &*=<3.5ex>[F]{7} &1 &3 &5 &*=<3.5ex>[F]{9} &6 &*=<3.5ex>[F]{10} &8
}}
$$ In the first row we box the left-to-right maxima to the right of
$1$ that are not right-to-left minima. Here, those are $7$ and $9$. In
the second row we insert a new largest letter, $10$, immediately to
the left of $1$ and box it. Finally, in the third row, we cyclically
shift the sequence of boxed letter one step to the left, thus
obtaining $\pi'$. Let us call this map $\beta$.

The induced map $\CK$, between $231$- and $321$-avoiding permutations
is then formally defined by
$$
\CK(\emptyword) = \emptyword;\quad
\CK(\alpha(\sigma)) = \beta(\CK(\sigma));\quad
\CK(\sigma\oplus\tau) = \CK(\tau)\oplus\CK(\sigma).
$$ 
As an example, consider the permutation $5213476$ in $\sym_6(231)$. Decompose
it using $\oplus$ and $\alpha$:
$$
5213476 
= 52134\oplus 21 
= \alpha(2134)\oplus\alpha(1) 
= \alpha(\alpha(1)\oplus 1\oplus 1)\oplus\alpha(1).
$$
Reverse the order of summands and change each $\alpha$ to $\beta$:
$$
\beta(1) \oplus \beta(1\oplus 1\oplus \beta(1)) 
= 21\oplus\beta(1243)
= 21\oplus 41253 
= 2163475.
$$
In conclusion, $\CK(5213476) = 2163475$.


\section{Proof of Theorem~\ref{relations}}
\label{sec:relations}

In the following five subsections we prove
Theorem~\ref{relations}---one subsection for each equality in the
theorem.

\subsection{Simion-Schmidt versus \CK}


We prepare for this proof by characterizing the Simion-Schmidt
bijection in terms of left-to-right minima. (That characterization can
be said to be implicit in \cite{SS85}.) We also characterize the
bijection \CK\ in terms of right-to-left minima.

\begin{definition}
  For a permutation $\pi = a_1a_2\dots a_n$ of length $n$, define 
  $$\LMIN(\pi) = \big\{\,(i,a_i) 
  \;\big\vert\;\text{$a_i$ is a left-to-right minima in $\pi$} \,\big\}
  $$ as the set of positions of left-to-right minima together with
  their values. Also, define $\sym_n/\LMIN$ as the set of equivalence
  classes with respect to the equivalence induced by $\LMIN$: that is,
  $\pi$ is equivalent to $\tau$ if\/ $\LMIN(\pi) =
  \LMIN(\tau)$. Similarly, define $\RMIN$, $\LMAX$ and $\RMAX$.
\end{definition}

The cardinality of $\LMIN(\sym_n)=\{\LMIN(\pi) \mid\pi\in\sym_n \}$ is
easily seen to be $C_n$, the $n$-th Catalan number. The following
lemma strengthens that observation:

\begin{lemma}\label{lemma:lmin}
  Each equivalence class in $\sym_n/\LMIN$ contains exactly one
  permutation that avoids $123$ and one that avoids $132$. In other
  words, both $\sym_n(123)$ and $\sym_n(132)$ are complete sets of
  representatives for $\sym_n/\LMIN$.
\end{lemma}

\begin{proof}
  Given a set $L$ in $\LMIN(\sym_n)$, this is how we construct the
  corresponding permutation $\tau=c_1c_2\dots c_n$ in $\sym_n(132)$:
  For $i$ from $1$ to $n$, if $(i,a)$ is in $L$ let $c_i = a$;
  otherwise, let $c_j$ be the smallest letter not used that is greater
  than all the letters used thus far.
  
  Given a set $L$ in $\LMIN(\sym_n)$, this is how we construct the
  corresponding permutation $\pi=a_1a_2\dots a_n$ in $\sym_n(123)$:
  For $i$ from $1$ to $n$, if $(i,c)$ is in $L$ let $a_i = c$;
  otherwise, let $a_j$ be the largest letter not used thus far.
  
  It is easy to see that filling in the letters in any other way than
  the two ways described will either change the sequence of
  left-to-right minima or result in an occurrence of $132$ or $123$.
\end{proof}

As an illustration of the preceding proof, with
$L=\{(1,6),(3,3),(4,2),(6,1)\}$ we get $67324158$ in $\sym_8(132)$ and
$68327154$ in $\sym_8(123)$.

Using Lemma~\ref{lemma:lmin} we can thus define a bijection between
$\sym_n(123)$ and $\sym_n(132)$ by letting $\pi$ correspond to $\sigma$
if $\LMIN(\pi) = \LMIN(\sigma)$. However, this map is not new---it
is the Simion-Schmidt bijection:

\begin{lemma}\label{lemma:SiS}
  For $\pi$ in $\sym_n(123)$ and $\sigma$ in $\sym_n(132)$, the
  following two statements are equivalent:
  \begin{enumerate}
  \item $\SiS(\pi) = \sigma$;
  \item $\LMIN(\pi) = \LMIN(\sigma)$.
  \end{enumerate}
\end{lemma}

Indeed, looking at the algorithm defining the Simion-Schmidt
bijection we see that the variable $x$ keeps track of the smallest
letter read thus far; lines 3 and 4 express that left-to-right minima
are left unchanged; and line 6 assign $c_j$ to be the smallest letter
not used that is greater than all the letters used thus far (as
described above).

Here is a characterization of $\CK$ in terms of $\RMIN$: 

\begin{lemma}\label{lemma:CK}
  For $\pi$ in $\sym_n(231)$ and $\sigma$ in $\sym_n(321)$, the
  following two statements are equivalent:
  \begin{enumerate}
  \item $\CK(\pi)=\sigma$;
  \item $(n+1-i,a)\in \RMIN(\pi) \iff (n+1-a,i)\in \RMIN(\sigma)$.
  \end{enumerate}
\end{lemma}

\begin{proof}
  That the latter statement characterizes a bijection from
  $\sym_n(231)$ to $\sym_n(321)$ follows from Lemma~\ref{lemma:lmin},
  so all we need to show is that $\CK$ is that bijection.

  We use induction on $n$, the length of the permutation. The case
  $n=1$ is obvious: the right-to-left minimum $(1,1)$ goes to the
  right-to-left minimum $(1,1)$. For the induction step we distinguish
  two cases: $\pi$ is decomposable and $\pi$ is indecomposable (see
  Section~\ref{sec:CK} for definitions).

  Suppose that $\pi$ is indecomposable, and hence $\pi =
  \alpha(\sigma)$ for some $\sigma$. By definition,
  $\CK(\pi)=\beta(\CK(\sigma))$. The claim follows from the
  following equivalences:
  \begin{align*}
    (n+1-i,a)\in \RMIN(\pi)
    &\iff (n-i,a)\in \RMIN(\sigma) \\
    &\iff (n-a,i)\in \RMIN \CK(\sigma) \\
    &\iff (n+1-a,i)\in \RMIN \CK(\pi).
  \end{align*}
  Here, the first equivalence is immediate from the definition of
  $\alpha$---recall that all $\alpha$ does is to insert a new largest
  letter in front of $\sigma$. The second equivalence holds by
  induction. The third equivalence follows from the definition of
  $\beta$: the cyclic shift involves no right-to-left minima and the
  space for the new letter is created to the left of the letter
  $1$; therefore, $1$ is added to the indices of right-to-left minima.

  Suppose $\pi=\tau\oplus\rho$ is decomposable, and let $k=|\tau|$
  and $\ell=|\rho|$. By definition,
  $\CK(\pi)=\CK(\rho)\oplus\CK(\tau)$. We have
  \begin{align*}
    &(n+1-i,a)\in\RMIN(\pi) \\
    &\iff (n+1-i,a)\in\RMIN(\tau)\text{ or }(\ell+1-i,a-k)\in\RMIN(\rho)
    &&\text{def' of $\oplus$} \\
    &\iff (k+1-(i-\ell),a)\in\RMIN(\tau)\text{ or }(\ell+1-i,a-k)\in\RMIN(\rho)
    &&\text{$n=k+\ell$} \\
    &\iff (k+1-a,i-\ell)\in\RMIN\CK(\tau)\text{ or }
    (\ell+1-(a-k),i)\in\RMIN\CK(\rho)
    &&\text{induction} \\
    &\iff (k+1-a,i-\ell)\in\RMIN\CK(\tau)\text{ or }(n+1-a,i)\in\RMIN\CK(\rho)
    &&\text{$n=k+\ell$} \\
    &\iff (n+1-a,i)\in\RMIN(\CK(\rho)\oplus\CK(\tau))
    &&\text{def' of $\oplus$}
  \end{align*}
  from which the claim follows.
\end{proof}

We now turn to the proof of the first identity of
Theorem~\ref{relations}. It is equivalent to
$$
\inverse\,\circ \,
\SiS\,\circ\,\reverse\,\circ\
\CK\,\circ\,\reverse\ =\ \identity.
$$
With all the preparation we have done, this is easy to prove:
\begin{align*}
  (i,a)\in\LMIN(\pi)
  &\iff (n+1-i,a)\in\RMIN.\reverse(\pi) \\
  &\iff (n+1-a,i)\in\RMIN.\CK.\reverse(\pi) \\
  &\iff (a,i)\in\LMIN.\reverse.\CK.\reverse(\pi) \\
  &\iff (a,i)\in\LMIN.\SiS.\reverse.\CK.\reverse(\pi) \\
  &\iff (i,a)\in\LMIN.\inverse.\SiS.\reverse.\CK.\reverse(\pi).
\end{align*}


\subsection{Simion-Schmidt versus Krattenthaler}

In Lemma~\ref{lemma:SiS} we characterized the Simion-Schmidt
bijection. We shall do the same for Krattenthaler's bijection. We
start by looking at the standard bijection from 132-avoiding
permutations to Dyck paths (as defined in Section~\ref{subsec:Knuth}).

\newcommand{\U}{\mathfrak{u}} 
\newcommand{\D}{\mathfrak{d}} 

Let $P$ be a Dyck path of length $2n$; index its up- and down-steps
$1$ through $n$. For instance,
$$P = u_1u_2u_3d_1d_2u_4u_5d_3d_4u_6d_5d_6.
$$ 
A \emph{peak} in a Dyck path is an up-step directly followed by a
down-step. Define
$$\PEAK(P) = \{\, (i,j) \mid u_id_j \text{ is a peak in }P \,\}.
$$
For instance, with $P$ as before, we have $\PEAK(P)=\{(3,1),(5,3),(6,5)\}$. 

\begin{lemma}\label{lemma:standard}
  Let $f$ be the standard bijection from $\sym_n(132)$ to $\dyck_n$.
  For $\pi$ in $\sym_n(132)$ and $P$ in $\dyck_n$, the following two
  statements are equivalent:
  \begin{enumerate}
    \item $f(\pi) = P$;
    \item $(i,n+1-a)\in\LMIN(\pi) \iff (a,i)\in\PEAK(P)$.
  \end{enumerate}
\end{lemma}

\begin{proof}
  Clearly, knowing $\PEAK(P)$ is equivalent to knowing the path
  $P$. Thus the second statement determines a bijection (by
  Lemma~\ref{lemma:lmin}). It remains to show that the first statement
  implies the second.

  We shall use induction on the length of the permutation. Assume that
  $\pi^r$ is indecomposable (with respect to $\oplus$). It is easy to
  see that $\pi$ ends with its largest letter. Hence, $\pi=\tau n$ for
  some $\tau$ in $\sym_{n-1}(132)$. Let $Q=f(\tau)$. Then
  $P=f(\pi)=uf(\tau)d=uQd$ and
  \begin{align*}
    (i,n+1-a)\in\LMIN(\pi)
    & \iff (i,n+1-a)=(i,n-(a-1))\in\LMIN(\tau) \\
    & \iff (a-1,i)\in\PEAK(Q) \\
    & \iff (a,i)\in\PEAK(P).
  \end{align*}

  Assume that $\pi^r$ is decomposable, so $\pi^r = \rho^r\oplus\tau^r$
  for some $\tau$ in $\sym_{k}(132)$ and $\rho$ in $\sym_{\ell}(132)$
  with both $k$ and $\ell$ positive and $k+\ell=n$.  Let
  $Q=f(\tau)$ and $R=f(\rho)$. Then $P=f(\pi)=f(\tau)f(\rho)=QR$ and
  \begin{align*}
    & (i,n+1-a)\in\LMIN(\pi) \\
    & \iff (i,k+1-a)\in\LMIN(\tau)\text{ or } (i-k,n+1-a)\in\LMIN(\rho)  \\
    & \iff (i,k+1-a)\in\LMIN(\tau)\text{ or } (i-k,\ell+1-(a-k))\in\LMIN(\rho)  \\
    & \iff (a,i)\in\PEAK(Q)\text{ or } (a-k,i-k)\in\PEAK(R) \\
    & \iff (a,i)\in\PEAK(P),
  \end{align*}
  which completes the proof.
\end{proof}

\begin{lemma}\label{lemma:g}
  Let $K$ be Krattenthaler's bijection from $\sym_n(123)$ to $\dyck_n$
  as described in Section~\ref{subsec:Krattan}.  For $\pi$ in
  $\sym_n(123)$ and $P$ in $\dyck_n$, the following two statements are
  equivalent:
  \begin{enumerate}
    \item $K(\pi) = P$;
    \item $(i,n+1-a)\in\RMAX(\pi) \iff (a,i)\in\PEAK(P)$.
  \end{enumerate}
\end{lemma}

\begin{proof}
  This is an easy consequence of $K$'s definition.
\end{proof}

Putting Lemma~\ref{lemma:standard} and Lemma~\ref{lemma:g} together we get
the desired characterization of Krattenthaler's bijection.

\begin{lemma}\label{lemma:Krattan}
  For $\pi$ in $\sym_n(123)$ and $\sigma$ in $\sym_n(132)$, the
  following two statements are equivalent:
  \begin{enumerate}
    \item $\Krattan(\pi) = \sigma$;
    \item $(n+1-i,a)\in\RMAX(\pi) \iff (n+1-a,i)\in\LMIN(\sigma)$.
  \end{enumerate}
\end{lemma}

Having established this characterization, the rest is easy. From the
sequence of equivalences
\begin{align*}
  &(i,a)\in\LMIN(\pi)  \\
  &\iff (n+1-i,a)\in\RMIN.\reverse(\pi) \\
  &\iff (a,n+1-i)\in\LMAX.\inverse.\reverse(\pi)\\
  &\iff (n+1-a,n+1-i)\in\RMAX.\reverse.\inverse.\reverse(\pi) \\
  &\iff (i,a)\in\LMIN.\Krattan.\reverse.\inverse.\reverse(\pi) \\
  &\iff (i,a)\in\LMIN.\SiS^{-1}.\Krattan.\reverse.\inverse.\reverse(\pi)
\end{align*}
it follows that
$$
\SiS^{-1}\circ\,\Krattan\,\circ\,\reverse\,\circ\,\inverse\circ\,\reverse 
\ =\ \identity
$$ 
as desired.

\subsection{Knuth-Richards versus \CK}

Consider the Dyck path $P=uudduududuuddd$ of semi-length $n=7$. Let us
index its up- and down-steps 1 through $n$:
$$P = u_1u_2d_1d_2u_3u_4d_3u_5d_4u_6u_7d_5d_6d_7.
$$ From this path we shall construct a permutation $\pi=a_1a_2\dots
a_n$. Scan $P$'s down-steps from left to right: if $d_i$ is preceded
by an up-step $u_j$, then let $a_{n+1-j}=i$; otherwise, let $j$ be the
largest value for which $a_j$ is unset, and let $a_j = i$. Like this:
\begin{enumerate}
  \item $d_1$ is preceded by the up-step $u_2$; let $a_{8-2}=a_6=1$.
  \item $d_2$ is preceded by a down-step; let $j=7$ and $a_7=2$.
  \item $d_3$ is preceded by the up-step $u_4$; let $a_{8-4}=a_4=3$.
  \item $d_4$ is preceded by the up-step $u_5$; let $a_{8-5}=a_3=4$.
  \item $d_5$ is preceded by the up-step $u_7$; let $a_{8-7}=a_1=5$.
  \item $d_6$ is preceded by a down-step; let $j=5$ and $a_5=6$.
  \item $d_7$ is preceded by a down-step; let $j=2$ and $a_2=7$.
\end{enumerate}
The resulting permutation is $\pi=5743612$. What we have just
described is the algorithm defining Richard's bijection. Lines 3, 4 and
5 of that algorithm covers the case when $d_i$ is preceded by an
up-step; lines 6, 7 and 8 the case when $d_i$ is preceded by a
down-step.

Plainly, if $d_i$ is preceded by an up-step $u_i$ then $u_id_j$ is a
peak in $P$. Moreover, $a_{n+1-j}=i$ is a left-to-right minimum in the
corresponding permutation. To be precise we have the following lemma.

\begin{lemma}
  Let $\Ri$ be Richards' bijection from $\sym_n(123)$ to $\dyck_n$ as
  described in Section~\ref{subsec:KRi}. For $\pi$ in $\sym_n(123)$
  and $P$ in $\dyck_n$, the following two statements are equivalent:
  \begin{enumerate}
    \item $\Ri(P) = \pi$;
    \item $(n+1-i,a)\in\LMIN(\pi) \iff (i,a)\in\PEAK(P)$.
  \end{enumerate} 
\end{lemma}

Using Lemma~\ref{lemma:standard} we get a characterization of the
Knuth-Richards bijection:

\begin{lemma}\label{lemma:KRi}
  For $\pi$ in $\sym_n(132)$ and $\sigma$ in $\sym_n(123)$, the
  following two statements are equivalent:
  \begin{enumerate}
    \item $\KRi(\pi) = \sigma$;
    \item $(i,a)\in\LMIN(\pi) \iff (a,i)\in\LMIN(\sigma)$.
  \end{enumerate}
\end{lemma}


The rest is easy. We have
\begin{align*}
(i, a) \in \LMIN(\pi)
& \iff (a,i) \in\LMIN.\KRi^{-1}(\pi)     \\
& \iff (n+1-a,i) \in\RMIN.\reverse.\KRi^{-1}(\pi) \\
& \iff (n+1-i,a) \in\RMIN.\CK.\reverse.\KRi^{-1}(\pi) \\
& \iff (i,a) \in\LMIN.\reverse.\CK.\reverse.\KRi^{-1}(\pi)
\end{align*}
and hence
$$\reverse\,\circ\,\CK\,\circ\,\reverse\,\circ\,\KRi^{-1}\, =\ \identity.
$$

\subsection{Simion-Schmidt versus Mansour-Deng-Du}

We will show that
$$
\MDD\,=\,\reverse\,\circ\,\SiS\,\circ\,\reverse.
$$ 
Due to Lemma~\ref{lemma:SiS} it suffices to prove this lemma:

\begin{lemma}
  For $\pi$ in $\sym_n(132)$ and $\sigma$ in $\sym_n(123)$, the
  following two statements are equivalent:
  \begin{enumerate}
  \item $\MDD(\pi) = \pi'$;
  \item $\RMIN(\pi) = \RMIN(\pi')$.
  \end{enumerate}
\end{lemma}


\begin{proof}

Assume that $\pi$ and $\pi'$ are as above. According to the proofs of
Corollaries~\cite[Cor. 4]{MDD} and~\cite[Cor. 16]{MDD}, the positions
of right-to-left minima in $\pi$ and $\pi'$ are the same and, in
particular, $\rmin(\pi)=\rmin(\pi')$. Thus we only need to prove that
right-to-left minima are preserved in value under the Mansour-Deng-Du
bijection. Equivalently, we need to prove that
non-right-to-left-minima ($\Nrmin$) are preserved in value.

One can see that a letter $a$ is an $\Nrmin$ in $\pi$ if and only if
the reduced decomposition of $\pi$ contains a run of simple
transpositions $(s_{a-1}\dots)$. In particular, $a=1$ is always an
$\Nrmin$. Thus $\pi=(12\dots n)\sigma_1\dots\sigma_k$ and
$\pi'=(12\dots n)\sigma'_1\dots\sigma'_k$ for $k=n-\rmin(\pi)$. That
is, $\pi$ and $\pi'$ have the same number of runs of simple
transpositions in the reduced decompositions and it remains to show
that the first letter of $\sigma_j$ equals the first letter of
$\sigma'_j$ whenever $1\leq j\leq k$.

Let $P$ be the intermediate Dyck path and consider its $(x+y)$- and
$(x-y)$-labellings of $P$. Note that cells touching the $x$-axis
receive the same labels under both labellings. From this, and the way
that the zigzag and trapezoidal decompositions are constructed, it
immediately follows that $\sigma_k$ and $\sigma'_k$ begin with the
same letter, namely, the label $C$ of the rightmost cell.

We now proceed by induction on the number of essential cells. If there
are no essential cells, then the statement is true. Suppose we have
$k>0$ essential cells. Remove the rightmost zigzag strip to get a Dyck
path $P'$. Note that $|P'| = |P|-2$ and that $P'$ has $k-1$ essential
cells. Clearly, the permutation corresponding to the $(x+y)$-labeling
of $P'$ is $\tau=(12\dots (n-1))\sigma_1\dots\sigma_{k-1}$. Let the
permutation corresponding to the $(x-y)$-labeling of $P'$ be
$\tau''=(12\dots (n-1))\sigma_1''\dots\sigma_{k-1}''$. From the
properties of the $(x-y)$-labeling, a cell labeled $Q\neq C$ is the
rightmost cell of a trapezoidal strip in $P$ if and only if $Q$ is the
rightmost cell of a trapezoidal strip in $P'$. This means that
$\sigma_i'$ and $\sigma_i''$ begin with the same letter for $1\leq
i\leq k-1$. The desired result now follows from the induction
hypothesis applied to $P'$, $\tau$ and $\tau''$.
\end{proof}

\subsection{Reifegerste versus Knuth-Rotem}

We will show that 
\begin{equation}\label{eq:Astrid_KRo}
\Astrid\,=\,
\inverse\,\circ\,\KRo\,\circ\,\inverse.
\end{equation}
Reifegerste's bijection is defined by values and positions of
excedances. Suppose that $\pi=a_1a_2\dots a_n$ is a 321-avoiding
permutation and that $\pi'$ is the image of $\pi$ under Reifegerste's
bijection. The only permutation of length $n$ without excedances is
$12\dots n$. As is easy to see, that permutation is fixed by both
sides of identity~\eqref{eq:Astrid_KRo}. So we can assume that $\pi$ has
at least one excedance. Say that has an excedance $a$ at position $i$,
denoted $(i,a)$. Also, let $(j,b)$ be the excedance closest to the
left of $(i,a)$. If no such excedance exists we define $j=0$. Note
that $b<a$, since otherwise an occurrence of $321$ would be formed.

Consider the Ferrer's diagram corresponding to $\pi$ in the definition
of Reifegerste's bijection (shaded in the figure in
Subsection~\ref{subsec:astrid} on page~\pageref{subsec:astrid}).  The
point $(j+1,n+2-a)$ is a corner in that diagram. It is sent to the
left-to-right minimum $(j+1,n+2-a)$ in $\pi'$. Thus, we need to prove
that an excedance $(i,a)$ corresponds to a left-to-right minimum
$(j+1,n+2-a)$ under the right hand side of
identity~\eqref{eq:Astrid_KRo}.

Because $a>i$, $b>j$ and $b<a$ it follows that both $i$ and $j$ are
{\em strict non-excedances} in $\pi^{-1}$ with the property that there
are no other strict non-excedances between $i$ and $j$. (A strict
non-excedance between $i$ and $j$ would either bring an occurrence of
321 in $\pi$ or in $\pi^{-1}$, or an occurrence of an excedance
between $a$ and $b$ in $\pi$.) Thus the ballot sequence $\beta$,
obtained when applying Knuth-Rotem's bijection to $\pi^{-1}$ will have
the letter $j+1$ in positions $b=a_j,a_{j+1},\dots,a_{i-1}$ and the
letter $i+1$ in position $a=a_i$. Let $P$ be the Dyck path
corresponding to $\beta$. It remains to show that $f^{-1}$, the
inverse of the standard bijection, sends $P$ to a permutation having
the letter $j+1$ in position $n+2-a$. After applying inverse we would
then have the letter $n+2-a$ in position $j+1$, the same outcome as
when applying Reifegerste's bijection.

\begin{figure}
\begin{center}
\vspace{20pt}

\setlength{\unitlength}{4mm}

\begin{picture}(25,6)
\put(5,0.5){

\path(0,5)(0,0)(5,0) \path(-0.2,4.7)(0,5)(0.2,4.7)
\path(4.7,0.2)(5,0)(4.7,-0.2)

\path(10,6)(10,0)(16,0) \path(9.8,5.7)(10,6)(10.2,5.7)
\path(15.7,0.2)(16,0)(15.7,-0.2)

\dashline{0.2}(0,0)(4.5,4.5)
\dashline{0.2}(10,0)(15.5,5.5)
\dashline{0.2}(11,0)(15.5,4.5)

\put(9.2,1.8){\small $i$} \put(9.2,0.8){\small $j$}
\path(9.9,2)(10.1,2)\path(9.9,1)(10.1,1) 
\put(13.8,-0.7){\small $a$}

\path(10,0)(11,0)(12,0)(12,1)(13,1)(14,1)(14,2)(15,2)(15,3)(15,4)(15,5)
\put(10,0){\p}\put(11,0){\p}\put(12,0){\p}
\put(12,1){\p}\put(13,1){\p}\put(14,1){\p}\put(14,2){\p}
\put(15,2){\p}\put(15,3){\p}\put(15,4){\p}\put(15,5){\p}

\path(-0.1,2)(0.1,2)  \linethickness{4mm}
\path(0,0)(2,0)(2,2)(4,2)(4,4) \put(-2.1,1.8){\small
$i=1$}\put(-2.1,-0.2){\small $j=0$} \put(1.05,-0.7){\small $a=2$}

\put(0,0){\p}\put(2,0){\p} \put(2,2){\p} \put(4,2){\p} \put(4,4){\p}

\put(2.5,-2){A} \put(12.5,-2){B}

}
\end{picture}
\ \vspace{5pt}

\smallskip 
\caption{}\label{proof001}
\end{center}
\end{figure}

For the remainder of this proof we use induction on $n$, the length of
the permutation $\pi$. The smallest permutation that have an excedance
is $\pi=21$. In this case, $i=1$, $j=0$ and $a=2$. See
Figure~\ref{proof001}A. After rotating that diagram counter clockwise
by $3\pi/4$ radians we read the Dyck path $udud$. The inverse of the
standard bijection, $f^{-1}$, sends $udud$ to the permutation
$\pi'=21$. It has the letter $j+1=1$ in position $n+2-a=2$ as desired.

Assume that $n>2$. Let $D$ be the diagram constructed from the
ballot-sequence corresponding to $\pi$. Let $P$ be the Dyck path we
read from $D$. Let $(r,s)$ be the coordinate in $D$ corresponding to
the first return to the $x$-axis in $P$. In particular, $r=s$.
Consider the following four cases.

{\it Case 1,} $s=0$. This case is sketched in
Figure~\ref{proof001}B. Note that $a>i+1$, and we can remove the first
step and the last step of the path, thus shifting the path one step to
the left. This corresponds to adjoining $n$ to the right side of the
permutation obtained inductively from the reduced Dyck path. Moreover,
the number of steps is now $2(n-1)$, $i$ and $j$ are unchanged, while
$a$ becomes $a-1$. From the induction hypothesis it follows that $j+1$
will be in position $(n-1)+2-(a-1)=n+2-a$, as desired.

{\it Case 2,} $n>s>i$. From the definition of the standard bijection
we see that the permutation $f^{-1}(P)$ is of the form $\pi'=\sigma
n\tau$ where each letter of $\sigma$ is larger than any letter of
$\tau$; the largest letter, $n$, is in position $n-s+1$; part
$\sigma$ is obtained inductively from the portion of the path above
the line $y=s$; and $\tau$ is obtained inductively from the portion of
the path below the line $y=s$. Applying the induction to $\tau$ we see
that the letter $j+1$ is in position $(s-1)+2-a=s+1-a$. Thus, in
$\pi'$, the letter $j+1$ is in position $(s+1-a)+(n-s+1)=n+2-a$.

{\it Case 3,} $0<s<i$. From the path we remove the part below the line
$y=s$, we remove the first and the last steps of the remaining path,
and we shift the obtained path so that it goes from $(n-s-1,n-s-1)$ to
$(0,0)$. The resulting path we call $P'$. It is responsible for
building the standardization $\sigma'$ of $\sigma$ in the outcome
permutation $\pi'=\sigma n\tau$. (Here $\sigma'$ is obtained from
$\sigma$ by adding $|\tau|$ to each of its letters.) Note that $i$,
$j$ and $a$ in $P$ become $i-s$, $j-s$ and $a-s-1$ in $P'$,
respectively. From the induction hypothesis applied to $P'$ it follows
that $(j-s)+1$ will be in position $(n-s-1)+2-(a-s-1)=n+2-a$ in
$\sigma'$. But $(j-s)+1$ in $\sigma'$ is $(j-s)+1+s=j+1$ in
$\sigma$. Thus, in $\pi'$, the letter $j+1$ will be in position
$n+2-a$.

{\it Case 4,} $s=i$. In this case $i+1=a$ and we need to prove that
the letter $j+1$ is in position $n+2-a=n+1-i$. Remove from $P$ the
first and last steps of the subpath from $(n,n)$ to $(i,i)$.  Shift
the path from $(n,n-1)$ to $(i+1,i)$ one step to the left. Let $P'$
denote the obtained path from $(n-1,n-1)$ to $(0,0)$. Now we apply
$f^{-1}$ to $P$ and get a 132-avoiding permutation of length $(n-1)$
of the form $\sigma\,(n-1)\,\tau$, where $|\tau|\geq i$. The first
return in $P'$ to the line $y=x$ will be either at $(i,i)$ or above
the line $y=i$. By the induction hypothesis, in $\sigma(n-1)\tau$, the
letter $j+1$ is in position $(n-1)+2-a=n-i$. This position is located
to the right of $n-1$. Looking at the definition of the standard
bijection we see that $f^{-1}(P)$ is $\sigma(n-1)\tau$ with the
largest letter, $n$, inserted immediately to the left of the letter
$j+1$, causing it to appear in position $n-i+1 = n+2-a$.

\section{Proof of Theorem~\ref{extensions}}\label{sec:statistics}

In this section we prove the equidistributions presented in
Theorem~\ref{extensions}. In light of Theorem~\ref{relations} it
suffices to consider the following five bijections:
$$\SiS,\, \Astrid,\, \West,\, \Knuth \text{ and } \ED.
$$ For brevity we shall in this section write
$\mathrm{stat}_1\simeq\mathrm{stat}_2$ when, for all $\pi$ in the
relevant domain, we have $\mathrm{stat}_1(\pi) =
\mathrm{stat}_2(\psi(\pi))$, where $\psi$ is the bijection under
consideration.

\subsection{Simion-Schmidt's bijection} 

Let $\pi$ be a 123-avoiding permutation of length $n$ and let $\pi'$
be the image of $\pi$ under the Simion-Schmidt bijection.  By
Lemma~\ref{lemma:SiS} we have $\LMIN(\pi)=\LMIN(\pi')$. In particular,
$\lmin \simeq \lmin$. A three letter segment $abc$ is a valley
precisely when $a$ and $b$ are left-to-right minima but $c$ is
not. Thus $\valley\simeq\valley$. Similarly, we have $\ldr\simeq\ldr$
and $\head\simeq\head$. Indeed, $\ldr$ is determined by the position
of the first non-left-to-right minimum and $\head$ is the first
left-to-right minimum.

$\bullet\ \compR \simeq \compR$: Suppose that $\pi=AB$ where $B$
is the rightmost reverse component. Then $B=mC$ where $m$ is a
left-to-right minimum in $\pi$. It follows that $\pi'=A'B'$ in which
$A'$ and $B'$ are the images, under Simion-Schmidt's bijection, of $A$
and $B$, respectively. Because the Simion-Schmidt bijection is an
involution it also follows that $B'$ is the rightmost reverse
component in $\pi'$. By induction on the number of reverse components
we thus have $\compR \simeq \compR$.

$\bullet\ \slmaxC\simeq \lir$: The statistic $\slmaxC$ is the position of the
second left-to-right minimum (it is defined to be $n$ if there is only
one left-to-right minimum).  Suppose $\pi=a_1Aa_2B$ where $a_1$ and
$a_2$ are the two leftmost left-to-right minima (the case $\pi=1\dots
n$ is trivial). Note that each letter in $A$ is larger than $a_1$ and
$a_2$. We know that $\pi'=a_1A'a_2B'$ and $|A|=|A'|$. To avoid the
pattern 132 the segment $A'$ must be increasing. Thus $\slmaxC\simeq
\lir$.

$\bullet\ \slmaxIR\simeq\lirI$: The statistic $\slmaxIR$ is one less
than the minimal $i$ such that the letter $i$ is to the left of the
letter 1. Note that such an $i$ in $\pi$ must be a left-to-right
minimum. Suppose $\slmaxIR(\pi)=i-1$ and $i\leq n$ (the case
$\slmaxIR(\pi)=n$ is trivial). Then $i$ and $1$ are two consecutive
left-to-right minima in $\pi$. Consequently, $i$ and $1$ are two
consecutive left-to-right minima in $\pi'$. Thus the letters
$2,3,\dots,i-1$ must be to the right of 1 in $\pi'$. To avoid forming
an occurrence of 132, those letters must also be in increasing
order. Thus $\lirI(\pi')=i-1$.

$\bullet\ \ldrI\simeq\ldrI$: By definition $\ldrI(\pi)$ is the
largest $i$ such that $i, i-1, \dots , 1$ is a subword in
$\pi$. In particular, $i+1$, if it exists, is to the right of $i$ in
$\pi$. Suppose that $\ldrI(\pi)=i$. Clearly, the letters $i, i-1,
\dots, 1$ are consecutive left-to-right minima in $\pi$ and $i+1$ is a
non-left-to-right minimum. Thus $i, i-1, \dots, 1$ is a subword in
$\pi'$ and $i+1$ is to the right of $i$ in $\pi'$; hence
$\ldrI(\pi')=i$.

$\bullet\ \headIR\simeq\rmin$: Plainly, $\headIR(\pi)=n-i+1$ where $i$
is the position of 1 in $\pi$.  Let $\pi' = \sigma 1\tau$. The letter
1 is in the same position in $\pi'$ as it is in $\pi$ and so
$|\tau|=n-i$. To avoid 132 the letters of $\tau$ must be in
increasing order. The sequence of right-to-left minima in $\pi'$ is
thus simply $1\tau$ and therefore $\headIR(\pi)=\rmin(\pi')$.

$\bullet\ \valleyI\simeq\valleyI$: By definition, $\valleyI(\pi)$ is the number
of letters $i$ in $\pi$ such that $i$ is to the left of both $i-1$ and
$i+1$. Suppose $i$ in $\pi$ is one of those letters counted by
$\valleyI(\pi)$. To avoid 123 the letter $i$ must be a left-to-right
minimum. Thus $i-1$ is a left-to-right minimum too, but $i+1$ is
not. Since $\LMIN(\pi)=\LMIN(\pi')$ this observation translates to
$\pi'$. That is, in $\pi'$, the letters $i$ and $i-1$ are
left-to-right minima whereas $i+1$ is not. Thus $i+1$ is to the right
of $i$ in $\pi'$ and the letter $i$ in $\pi'$ is counted by
$\valleyI(\pi')$. Because the Simion-Schmidt bijection is an
involution it is easy to see that no $i$ which is not counted by
$\valleyI(\pi)$ contributes to $\valleyI(\pi')$.

$\bullet\ \rank\simeq\rank$: Suppose that $\rank(\pi)=k$ and $\pi$ has
$a$ in position $k+1$. We distinguish two cases based on whether $a$
is a left-to-right minimum.  If $a$ is a non-left-to-right-minimum,
then $k+1$ is the left-to-right minimum closest to $a$ in $\pi$ from
the left.  Since $\LMIN(\pi)=\LMIN(\pi')$ we have $\rank(\pi')=k$.  On
the other hand, if $a$ is a left-to-right minimum, then $a\leq k+1$
and $\pi'$ will have the same left-to-right minimum, $a$, in position
$k+1$. Thus $\rank(\pi')=k$ in this case as well.

\subsection{Reifegerste's bijection}

Let $\pi=a_1a_2\dots a_n$ be a 321-avoiding permutation of length $n$
and let $\pi'$ be the image of $\pi$ under Reifegerste's bijection.
That $\exc\simeq\des$ was proved in~\cite{R02}.

$\bullet\ \valley\simeq\valley$: A 321-avoiding permutation is a
shuffle of two increasing subwords. From this one can see that if
$a_ia_{i+1}a_{i+2}$ is a valley in $\pi$, then $a_i$ is an excedance
and $a_{i+1}$ is a non-excedance. Thus, in the permutation matrix
corresponding to $\pi$, there is an $E$-square in the $i$-th row but
no $E$-square in the $(i+1)$-th row. It follows that the Ferrer's
diagram has a corner in the $(i+1)$-th row. Thus
$a'_ia'_{i+1}a'_{i+2}$ is a valley in $\pi'=a'_1a'_2\dots a'_n$. (The
dot in row $i+1$, corresponding to the letter $a'_{i+1}$, will be to
the left of the dots in rows $i$ and $(i+2)$, corresponding to the
letters $a'_i$ and $a'_{i+2}$, respectively).

$\bullet\ \peakI\simeq\valleyI$: By definition, $\peakI(\pi)$ is the
number of letters $a$ in $\pi$ to the right of both $a-1$ and $a+1$;
and $\valleyI(\pi)$ is the number of letters $a$ in $\pi$ to the left
of both $a-1$ and $a+1$. If $i$ is counted by $\peakI$, then $i$ is a
non-excedance while $i+1$ is an excedance; otherwise an occurrence of
$321$ is formed. Thus we have an $E$-square in column $n-i$, but no
$E$-square in column $n-i+1$. Not also that there is a column $n-i+2$
(corresponding to the letter $i-1$). Thus column $n-i+1$ is not the
rightmost column of the matrix and it contains a corner of the
Ferrer's diagram. So $\peakI$ counts the columns that contain a corner
but no $E$-squares, excluding the rightmost column. In the
construction of $\pi'$ each corner with the properties described above
gives rise to a letter $i$ in $\pi'$ to the left of $i-1$ and
$i+1$. Indeed, such a corner has an $E$-square in the column
immediately to its left and no $E$-square in the column immediately to
its right; consequently, the points corresponding to $i-1$ and
$i+1$ in $\pi'$ will be below the point corresponding to $i$ in
$\pi'$. Thus $\peakI\simeq\valleyI$.

$\bullet\ \slmaxI\simeq\zeil$: The letter between the two leftmost
left-to-right maxima in a 321-avoiding permutation is an initial
segment of $123\dots$. Hence $\slmaxI$ is the length of the maximal
initial segment of the form $23\dots i$. Moreover, each of the
letters counted by $\slmaxI$ is an excedance. Thus we have
$E$-squares in positions $(\ell,n-\ell)$ for $\ell=1,2,\dots,i-1$ and
no $E$-square in position $(k,n-k)$. Thus $\pi'$ is of the form
$\pi'=n(n-1)(n-2)\dots (n-k+2)A(n-k)B(n-k+1)C$, where $A$, $B$, and
$C$ are some words. Clearly, $\zeil(\pi')=k=\slmaxI(\pi)$.

$\bullet\ \headI\simeq\ldr$: The statistic $\headI(\pi)$ is the
position of 1 in $\pi$. Say this position is $i$. Then row $i$ of the
permutation matrix corresponding to $\pi$ is the topmost row that does
not contain an $E$-square. Thus, in $\pi'$, the first $i$ letters
will be in decreasing order, while in position $i$ there will be an
ascent (unless $i=n$). Thus $\ldr(\pi')=i$.

$\bullet\ \slmaxRC\simeq\rdr$: Let the right-to-left minima of $\pi$,
read from right to left, be $r_1$, $r_2$, $\dots$, $r_k$.  One can see
that $\slmaxRC$ is one more than the number of letters between $r_1$
and $r_2$. (To make sure there is at least two left-to-right minima we
can assume that a 0 stays in front of $\pi$ when considering this
statistic.) To avoid the pattern 321 the letters between $r_1$ and
$r_2$ must be the largest letters in $\pi$, and thus all of them are
excedances. Thus, in the permutation matrix corresponding to $\pi$,
there will be $E$-squares in positions $(i,n-i)$ for
$i=n-1,n-2,\dots,n-(\slmaxRC(\pi)-1)$ and there will be no $E$-squares
in the ($n-\slmaxRC(\pi)$)-th row. In turn, this guarantees that the
132-avoiding permutation $\pi'$ ends with
$j$($\slmaxRC(\pi)-1$)($\slmaxRC(\pi)-2$)$\dots 1$, where $j$, if it
exists, is strictly larger than $\slmaxRC(\pi)$. To avoid 132 this
segment must be preceded by a letter smaller than $j$. Thus
$\rdr(\pi')=\slmaxRC(\pi)$.

$\bullet\ \rir\simeq\rmin$: Let $\rir(\pi)=i$. We have $i=n$ only if
$\pi=\pi'=12\dots n$, and then, trivially, $\rir(\pi) = \rmin(\pi') =
n$. Assume $i<n$. The letters in the rightmost increasing run are
non-excedances while the letter immediately to the left of that run
is an excedance. Thus the bottom-most $E$-square is in row $n-i$ in
the permutation matrix. The dots placed in rows $n-i+1$, $n-i+2,\dots,
n$ when creating $\pi'$ gives the sequence of right-to-left minima in
$\pi'$, and thus $\rmin(\pi')=i$.

$\bullet\ \lirI\simeq \lmax$: Let $\lirI(\pi)=i$. We have $i=n$ only if
$\pi=\pi'=12\dots n$, and then, trivially, $\lirI(\pi) = \lmax(\pi') =
n$. Assume $i<n$. By definition, $i$ is the largest positive integer
such that $1,2,\dots, i$ is a subword of $\pi$.  Since $\pi$ is
the shuffle of two increasing sequences, $i+1$ is the leftmost
excedance in $\pi$. The $E$-square corresponding to $i+1$ in the
permutation matrix is placed in column $n-i$, leaving $i$ columns to
the right of it. When constructing $\pi'$ each of those $i$ columns
get a dot, beginning at the position $(1,n-i+1)$ and going in the
South-East direction. Those dots give the sequence of left-to-right
maxima in $\pi'$, and thus $\lmax(\pi')=i$.

$\bullet\ \last\simeq\MldrI$: Note that $\ldrI(\pi)$ is the largest
$i$ such that $i, (i-1),\dots, 1$ is a subword in $\pi$. The case
$\pi=\pi'=12\dots n$ is trivial. Suppose $i<n$ and $\last(\pi)=i$. To
avoid the pattern 321, the letters $i+1$, $i+2,\dots,n$ must form a
subword of $\pi$, and clearly each of them is an excedance.
Therefore the $n-i$ leftmost columns of the permutation matrix
corresponding to $\pi$ contain $E$-squares but the $(n-i+1)$-th
column does not contain an $E$-square. From this, and the way $\pi'$
is constructed, it immediately follows that $(n-i+1),(n-i),\dots,1$ is
the longest subword of $\pi'$ of the sought type. Thus
$\MldrI(\pi')=i$.

\subsection{West's bijection}
Recall that West's bijection is induced by an isomorphism between
generating trees. The two isomorphic trees generate 123- and
132-avoiding permutations, respectively. However, for the purpose of
this proof we shall generate 321-avoiding permutations instead of
123-avoiding ones. This change is reflected in positions of active
sites: In the 123-avoiding case the active sites are all the positions
to the left of the leftmost ascent and the position in between the two
letters of the leftmost ascent. In the 321-avoiding case the active
sites are all the positions to the right of the rightmost ascent and
the position in between the two letters of the rightmost descent.

The active sites in a 132-avoiding permutation are the leftmost
position and every position immediately to the right of right-to-left
maxima.

All the proofs below are by induction on the length of the
permutation, with easily verifiable base cases. Let $\pi$ be a
321-avoiding permutation of length $n$. Let $\pi'$ be the image of
$\pi$ under the modified version of West's bijection as described
above.

$\bullet\ \peakI\simeq\valleyI$: By definition, $\peakI(\pi)$ is the
number of letters $a$ in $\pi$ to the right of both $a-1$ and $a+1$;
and $\valleyI(\pi')$ is the number of letters $a$ in $\pi'$ to the
left of both $a-1$ and $a+1$.

{\em Case 1}. Assume that $\pi$ ends with $n$. Then $\pi'$ begins with
$n$. Inserting $n+1$ to the right of $n$ in $\pi$ does not change
$\peakI$. Inserting $n+1$ to the left of $n$ in $\pi'$ does not change
$\valleyI$. Inserting $n+1$ in any other position in $\pi$ and $\pi'$
increases $\peakI$ and $\valleyI$ by 1.

{\em Case 2.} Assume that $\pi$ does not end with $n$. To avoid the
pattern 321, the letter $n$ must be the the last letter of the
rightmost descent. Thus $n$ is not the leftmost letter in $\pi'$, and
to avoid the pattern 132, the letter $n-1$ must be to the left of
$n$. It follows that inserting $n+1$ in an active site of $\pi$ and
$\pi'$, respectively, does not change $\peakI$ and
$\valleyI$, respectively.

$\bullet\ \exc\simeq\asc$: Inserting $n+1$ at the end of $\pi$ does
not change the number of excedances. Similarly, inserting $n+1$ at the
beginning of $\pi'$, the number of ascents is not changed. In all other
cases, the number of excedances in $\pi$ and the number of ascents in
$\pi'$ is increased by 1.

$\bullet\ \slmaxI\simeq\lirI$: We see that $\slmaxI(\pi)$ is one more
than the length of the maximum initial segment of the form $234\dots$
in $\pi$. Also, $\lirI(\pi')$ is the largest $i$ such that $12\dots i$ is
a subword of $\pi'$.

{\em Case 1.} Assume that $\pi=234\dots n1$. Using induction it is
easy to verify that $\pi'=12\dots n$. Inserting $n+1$ at the end of
$\pi$ does not change $\slmaxI$, and inserting $n+1$ at the beginning
of $\pi'$ does not change $\lirI$. On the other hand, inserting $n+1$
between $n$ and 1 in $\pi$ increases $\slmaxI$ by 1, and inserting
$n+1$ at the end of $\pi'$ increases $\lirI$ by 1.

{\em Case 2.} Assume that $\pi\neq 234\dots n1$, and thus $\pi'\neq
12\dots n$. Inserting $n+1$ will not change $\slmaxI$ in $\pi$ and it
will not change $\lirI$ in $\pi'$.

$\bullet\ \slmaxRC\simeq\comp$: By definition $\slmaxRC(\pi)=1$ if
the letter $n$ is not in position $n-1$, and $\slmaxRC(\pi)$ is one
more than the length of the maximal segment of the form $i(i+1)\dots
n$ if $n$ is in position $n-1$.

Inserting $n+1$ at the end of $\pi'$ (creating two active sites)
increases the number of components by 1. The corresponding action on
$\pi$ (creating two active sites) is inserting $n+1$ in position $n$;
that increases $\slmaxRC$ by 1.  Inserting $n+1$ in any other active
site of $\pi'$ will create an indecomposable permutation (n+1 will be
to the left of 1). The corresponding operation on $\pi$ places $n+1$
in a position different from $n$, and thus $\slmaxRC$ will be 1.

$\bullet\ \rir\simeq\rmax$: A proof is straightforward from the
location of active sites in $\pi$ and $\pi'$: the (number of) active
sites in $\pi$ and $\pi'$ essentially give the $\rir$- and
$\rmax$-statistics, respectively.

$\bullet\ \lirI\simeq\ldrI$: By definition $\lirI(\pi)$ is the largest
$i$ such that $12\dots i$ is a subword of $\pi$, and $\ldrI(\pi')$
is the largest $i$ such that $i(i-1)\dots 1$ is a subword of $\pi'$.

It is easy to see, by induction, that $\pi=12\dots n$ corresponds to
$\pi'=n(n-1)\dots 1$. Inserting $n+1$ at the end of $\pi$ and at the
beginning of $\pi'$ will increase $\lirI$ and $\ldrI$ by $1$,
respectively. For any other insertions and any other $\pi$ and $\pi'$,
the statistics $\lirI$ and $\ldrI$ will not change.

$\bullet\ \last\simeq\head$: Inserting $n+1$ in $\pi$ changes the last
letter only if the insertion is at the end. Similarly, inserting
$n+1$ in $\pi'$ changes the leftmost letter only if the insertion is
at the beginning.

\subsection{Knuth's bijection} 
Elizalde and Pak~\cite{EP04} proved that $\exc\simeq\exc$,
$\fp\simeq\fp$ and $\lis\simeq\Nrank$. Let $\pi$ be a $321$-avoiding
permutation of length $n$ and let $\pi'$ be the image of $\pi$ under
Knuth's bijection.

$\bullet\ \lir\simeq\lmax$: If $\lir(\pi)=i$ then the first row of the
{\em recording} tableau begins with $1,2,\dots, i$, and $i+1$ (if it
exists) is the leftmost element in the second row. Thus the statistic
$\lir$ translates to the statistic ``length of the rightmost slope''
(segment of down-steps) in the corresponding Dyck path. After
reflection, that statistic becomes ``length of the leftmost slope''
(segment of up-steps). The up-steps in the leftmost slope have
corresponding down-steps such that between these steps one has a
proper Dyck path. In particular, the down-step corresponding to the
leftmost up-step gives the first return to $x$-axis, giving the
position of $n$, the rightmost left-to-right maxima, in
$\pi'$. Proceeding recursively it is easy to see that in general the
up-steps in the leftmost slope read from left to right correspond to
the left-to-right maxima in $\pi'$ read from right to left. This gives
the desired result.

$\bullet\ \lirI\simeq\rmin$: The statistic $\lirI$ is the length of
the longest subword of the form $12\dots i$. So, if
$\lirI(\pi)=i$, the first row of the {\em insertion} tableau begins
with $1,2,\dots, i$, and $i+1$ (if it exists) is the leftmost element
in the second row. Thus the statistic $\lirI$ translates to the
statistic ``length of the leftmost slope'' in the corresponding Dyck
path.  After reflection, that statistic becomes ``length of the
rightmost slope''. Returns to $x$-axis in the Dyck path correspond to
reverse components in $\pi'$.  Consider the part $D'$ of the Dyck path
between the last return and next to last return to $x$-axis. The first
up-step of $D'$ corresponds to the rightmost down-step in the
rightmost slope, and it corresponds to the rightmost letter, say $a$,
in $\pi'$. The letter $a$ is the largest letter in the rightmost
reverse component of $\pi'$, and thus it is a right-to-left
minimum. Proceeding recursively we see that the next to last
down-step in the rightmost slope corresponds to the second
right-to-left minimum from the right, and so on. 

\subsection{Elizalde-Deutsch's bijection}

Elizalde and Deutsch~\cite{ED03} proved that $\fp\simeq\fp$.


\end{document}